\documentclass[final,onetabnum]{siamart171218}

\usepackage{lipsum}
\usepackage{amsfonts}
\usepackage{bm}   
\usepackage{amsmath}
\usepackage{graphicx}
\usepackage{dsfont}
\usepackage{epstopdf}
\usepackage{algorithmic}
\usepackage{booktabs}
\usepackage{enumerate}
\usepackage[toc]{appendix}
\usepackage{multirow}
\usepackage{soul}
\usepackage{caption}
\usepackage{subcaption}
\usepackage{scalerel,stackengine}
\newcommand{\reallywidehat}[1]{%
\savestack{\tmpbox}{\stretchto{%
\scaleto{%
\scalerel*[\widthof{\ensuremath{#1}}]                         {\kern-.6pt\bigwedge\kern-.6pt}%
{\rule[-\textheight/2]{1ex}{\textheight}}
}{\textheight}%
}{0.5ex}}%
\ensurestackMath{\stackon[1pt]{#1}{\tmpbox}}%
}

\ifpdf
  \DeclareGraphicsExtensions{.eps,.pdf,.png,.jpg}
\else
  \DeclareGraphicsExtensions{.eps}
\fi

\newtheorem{example}{Example}
\newtheorem{remark}[theorem]{Remark}
\def\RR{\mathbb{R}}

\def\BU{{\bm U}}
\def\BV{{\bm V}}
\def\BY{{\bm Y}}

\def\BL{{\bm L}_k}
\def\BPhi{{\bm \Phi}}
\def\BI{{\bm I}}
\def\BW{{\bm W}}
\def\bu{{\bm u}}
\def\by{{\bm y}}

\def\MPOD{{{\sc 2s-pod}}}
\def\MODEIM{{{\sc 2s-deim}}}
\def\MPDEIM{{{\sc 2s-pod-deim}}}
\def\podeim{{{\sc pod-deim}}}

\definecolor{greenf}{rgb}{.054, .5, .005}

\newcommand{\In}{\mathsf{\Phi}_{\bm f}}

\headers{POD-DEIM for semilinear matrix differential equations}{G. Kirsten and V. Simoncini}

\title{A matrix-oriented POD-DEIM algorithm applied to semilinear matrix differential equations  
\thanks{This version dated \today}}

\author{Gerhard Kirsten\thanks{Dipartimento di Matematica,
Universit\`a di Bologna, Piazza di Porta S. Donato, 5, I-40127 Bologna, Italy,
({\tt gerhard.kirsten2@unibo.it})}
\and
Valeria Simoncini\thanks{Dipartimento di Matematica and AM$^2$,
Universit\`a di Bologna, Piazza di Porta S. Donato, 5, I-40127 Bologna, Italy;
and IMATI-CNR, Pavia Italy ({\tt valeria.simoncini@unibo.it}).}}

\usepackage{amsopn}

\begin{document}

\maketitle

\begin{abstract}
We are interested in numerically approximating the solution $\BU (t)$
of the large dimensional semilinear matrix differential equation
$\dot{\BU}(t) = { \bm A}\BU (t) + \BU (t){ \bm B} + {\cal F}(\BU,t)$, 
 with appropriate starting and boundary
conditions, and $ t \in [0, T_f]$.
In the framework of the Proper Orthogonal Decomposition (POD) methodology and the
Discrete Empirical Interpolation Method (DEIM), we
derive a novel matrix-oriented reduction process leading to an effective,
structure aware low order approximation of the original problem.
The reduction of the nonlinear
term is also performed by means of a fully matricial interpolation using left and right
projections onto two distinct reduction spaces, giving rise to a new two-sided version of DEIM.
By maintaining a matrix-oriented reduction, we are able to employ first order exponential integrators
at negligible costs.
 Numerical experiments on benchmark problems illustrate the effectiveness of
the new setting. 
\end{abstract}


 
\section{Problem description}\label{sec:intro}
We are interested in numerically approximating the solution $\BU (t) \in {\cal S}$   
to the following semilinear matrix differential equation
 \begin{equation} \label{Realproblem}
\dot{\BU}(t) = {\color{black} \bm A}\BU (t) + \BU (t){\color{black} \bm B} + {\cal F}(\BU,t) , \quad \BU(0) = \BU_{0},
\end{equation}
where ${\color{black} \bm A} \in \mathbb{R}^{n_{x} \times n_{x}}, {\color{black} \bm B} \in \mathbb{R}^{n_{y} \times n_{y}}$,
and $ t \in [0, T_f] = {\cal T} \subset \RR$, equipped with appropriate boundary conditions.
The function {${\cal F} : {\cal S} \times {\cal T} \rightarrow \mathbb{R}^{n_x \times n_y}$} is a sufficiently 
regular nonlinear function that can be evaluated elementwise, and ${\cal S}$ is a functional space 
containing the sought after solution. 

The problem (\ref{Realproblem}) arises for instance in the discretization 
of two-dimensional partial differential equations of the form
\begin{equation} \label{nonlinpde}
u_{t} = \ell(u) + f(u,t) , \quad u = u(x,y, t) \quad \mbox{with} \,\,
 (x,y) \in \Omega \subset \mathbb{R}^2,\, t\in {\cal T},
\end{equation}
and given initial condition $u(x, y, 0) = u_{0}(x,y)$, for certain choices of the physical 
domain $\Omega$.
The differential operator $\ell$ is linear in $u$, typically a second order
operator in the space variables, while 
{$f: S \times {\cal T} \rightarrow \RR$} is a nonlinear function,
 where $S$ is an appropriate space with $u\in S$.
Time dependent equations of type (\ref{nonlinpde}) 
arise in biology and ecology, chemistry and physics,
where the interest is
in monitoring the time evolution of a complex phenomenon;
see, e.g.,
\cite{Mainietal.01},%
\cite{Malchowetal.08},%
\cite{Quarteroni.17},%
\cite{Tveitoetal.10}, and references therein.%

We develop a matrix-oriented \podeim\ order reduction strategy for the problem (\ref{Realproblem})
that leads to a semilinear {\it matrix} differential equation with the same
structure as (\ref{Realproblem}), but of significantly reduced dimension.
{More precisely, we determine an approximation to  $\BU(t)$ of the type 
\begin{equation}\label{eqn:Uapprox}
{\bm  V}_{\ell,U}{\bm Y}_k(t){\bm  W}_{r,U}^{\top} , \quad t \in [0, T_f] ,
\end{equation}
where ${\bm  V}_{\ell,U}\in \mathbb{R}^{n \times k_1}$ and 
${\bm  W}_{r,U}\in \mathbb{R}^{n \times k_2}$ are matrices to be determined, 
independent of time. Here $k_1, k_2 \ll n$ and we let $k=(k_1,k_2)$.
The function ${\bm Y}_k(t)$ is determined as the numerical solution to the following 
{\it reduced} semilinear matrix differential problem
 \begin{equation} \label{Realproblemsmall0}
 \begin{split}
\dot{{\bm Y}}_k(t) &= 
{\bm A}_{k}{\bm Y}_k(t) + {\bm Y}_k(t){\bm B}_{k} + 
                      \reallywidehat{{\cal F}_{k}({\bm Y}_k,t)} \\ 
{\bm Y}_k(0) &= {\bm Y}_k^{(0)} := {\bm  V}_{\ell,U}^{\top}\BU_{0}{\bm  W}_{r,U} ,
\end{split}
\end{equation}
with
${\bm A}_{k} = {\bm  V}_{\ell,U}^{\top}{\color{black} \bm A}{\bm  V}_{\ell,U}$, 
${\bm B}_{k} = {\bm  W}_{r,U}^{\top} {\color{black} \bm B}{\bm  W}_{r,U},$ 
and
$\reallywidehat{{\cal F}_{k}({\bm Y}_k,t)}$ is a matrix-oriented DEIM approximation 
to 
\begin{equation}\label{eqn:F}
{\cal F}_{k}({\bm Y}_k,t)=
{\bm  V}_{\ell,U}^{\top}{\cal F}({\bm  V}_{\ell,U}{\bm Y}_k{\bm  W}_{r,U}^{\top},t){\bm  W}_{r,U}.
\end{equation}

Standard procedures for (\ref{Realproblem}) employ a vector-oriented approach: 
 semi-discretization of \cref{nonlinpde} in space 
leads to the following system of  ordinary differential equations (ODEs) 
 \begin{equation} \label{vecode}
 \dot{\bu}(t) = {\bm L}\bu (t) + {\bm f}(\bu,t) , \quad \bu(0) = \bu_{0}.
 \end{equation}
For $t>0$, the vector $\bu (t)$ contains the representation coefficients of the sought after solution in the
chosen discrete space,
 ${\bm L} \in\RR^{N\times N}$ accounts for the discretization of the linear differential operator $\ell$
and {${\bm f}$} is 
evaluated componentwise at ${\bm u}(t)$.  

  The discretization of  (\ref{nonlinpde}) can directly lead to the form \cref{Realproblem}
whenever $\ell$ is
a second order differential operator with separable coefficients and a tensor basis is
used explicitly or implicitly for its
discretization, such as finite differences on structured grids, 
certain spectral methods, isogeometric analysis, etc.
In a tensor space discretization, for instance,
the discretization of
$\ell(u)$ yields the operator $\BU \mapsto {\bm  A} \BU + \BU {\bm  B}$, where
the matrices ${\color{black} \bm A}$ and ${\color{black} \bm B}$ approximate
the second order derivatives in the 
$x-$ and $y-$directions, respectively, using $n_x$ and $n_y$ discretization 
nodes\footnote{Here we display the discretized Laplace operator, more
general operators can be treated; see, e.g., \cite[Section 3]{Simoncini2017}.}.
%
%
%
We aim to show that by sticking to the matrix formulation of the problem throughout
the computation - including the reduced model - 
major advantages in terms of memory allocations, computational costs and
structure preservation can be obtained.

For ${\cal F}\equiv 0$, the equation (\ref{Realproblem}) simplifies to the differential
Sylvester equation, for which
reduction methods have shown to be competitive; see, e.g., 
\cite{Behr.Benner.Heiland.19} and references therein. 
Our approach generalizes to the semilinear case a matrix-oriented projection methodology successfully 
employed in this linear and quadratic contexts \cite{Kirsten2019}.
A challenging difference lies in the selection and construction of the
two  matrices ${\bm  V}_{\ell,U}, {\bm  W}_{r,U}$, so as to effectively handle
the nonlinear term in (\ref{Realproblem}). In addition, ${\cal F}$ itself
needs to be approximated for efficiency.

Order reduction of the vector problem (\ref{vecode})  is a well established procedure. Among various
methods, 
the Proper Orthogonal Decomposition (POD) methodology has been widely employed, as it mainly
relies on solution samples, rather than the a-priori generation of an appropriate basis
\cite{benner2015review},\cite{Benner.05},\cite{hinze2005},\cite{kunisch1999}.
Other approaches include
reduced basis methods, see, e.g., \cite{patera2007reduced}, and rational interpolation strategies 
\cite{ABG.20}; see, e.g., \cite{BCOW.17}
for an overview of the most common reduction strategies.
{The overall effectiveness of the POD procedure is largely influenced by the capability
of evaluating the nonlinear term within the reduced space, motivating a
considerable amount of work towards this estimation, including 
quadratic  bilinear approximation \cite{gu2011,kramer2019nonlinear, benner2015}
and trajectory piecewise-linear approximation \cite{white2003}. 
 Alternatively several approaches consider interpolating the nonlinear function,
such as 
missing point estimation \cite{astrid2008} and the best points interpolation 
method \cite{nguyen2008}. 
One very successful approach is the Discrete Empirical Interpolation Method 
(DEIM) \cite{chaturantabut2010nonlinear}, which is based on the Empirical Interpolation Method 
originally introduced in \cite{barrault2004}. 

We devise a matrix-oriented POD approach tailored
towards the construction of the
matrix reduced problem formulation (\ref{Realproblemsmall0}).
{An adaptive procedure is also developed to limit the
number of snapshots contributing to the generation of the approximation spaces.}
The reduction of the nonlinear term is then performed by means of
a fully matricial interpolation using left and right projections
onto two distinct reduction spaces, giving rise to a new two-sided version of DEIM. 

The idea of using left and right POD-type bases in a matrix-oriented setting is not new
in the general context of semilinear differential equations (see section~\ref{sec:others}). Nonetheless, 
after reduction these strategies resume the vector form of (\ref{Realproblemsmall0}) for integration
purposes, thus loosing the structural and computational benefits of the matrix formulation.
We claim that once (\ref{Realproblemsmall0}) is obtained, matrix-oriented integrators should be
employed. 
 In other words,
by combining matrix-oriented versions of POD, DEIM and ODE integrators, we are able to
carry the whole approximation with explicit reference to the two-dimensional computational domain.
%
As a result, a fast (offline) reduction phase where a significant decrease in the 
problem size is carried out, is
followed by a light (online) phase where the reduced ordinary differential matrix equation is
integrated over time with the preferred matrix-oriented method. 

Our construction focuses on the two-dimensional problem. The advantages of our methodology
become even more apparent in the three-dimensional (3D) case. A simplified version of our
framework in the 3D case is experimentally explored in the companion
manuscript \cite{Kirsten.21}, where
the application to systems of differential equations is also discussed. Here we deepen
the analysis of all the ingredients of this new methodology, and emphasize its advantages
over the vector-based approaches with a selection of numerical results. A more extensive
experimental evidence can be found in the previous version of this paper \cite{Kirsten.Simoncini.arxiv2020}.

The paper is organized as follows. 
In \cref{sec:POD_DEIM} we review the standard \podeim\ algorithm for systems of 
the form \cref{vecode}, whereas our new two-sided proper orthogonal decomposition is derived in \cref{sec:matpod}. In \cref{sec:others} we discuss the relation to other matrix-based interpolation strategies and in \cref{sec:dynamic} 
we present a dynamical procedure for selecting the snapshots. Section~\ref{sec:matdeim} is devoted to the crucial approximation of the
nonlinear function by the new two-sided discrete empirical interpolation method.
 The overall new procedure with the numerical treatment of the reduced differential problem
is summarized in \cref{extend}.
Numerical experiments are reported in
\cref{sec:exp} to illustrate the effectiveness of the proposed procedure.
Technical implementation details and computational costs are discussed in \cref{sec:alg}.

{\it Notation.} 
${\bm I}_n$ denotes the $n \times n$ identity matrix{; the subscript is omitted
whenever clear from the context}. For a matrix ${\bm A}$,
$\|{\bm A}\|$ denotes the matrix norm induced by the Euclidean vector norm, and 
$\|{\bm A}\|_F$ the Frobenius norm. Scalar quantities are denoted by lower case letters, while vectors (matrices)
are denoted by bold face lower (capital) case letters. 
As an exception, matrix quantities of interest in the {\it vector} POD-DEIM approximation 
 are denoted in sans-serif, instead of bold face, i.e., ${\sf M}$. 

All reported experiments were performed using MATLAB 9.6 (R2020b) 
(\cite{matlab2013}) on a MacBook Pro with 8-GB memory and a 2.3-GHz Intel core i5 processor.

\section{The standard POD method and DEIM in the vector framework}\label{sec:POD_DEIM}
We review the standard \podeim\ method and its application to the 
dynamical system \cref{vecode}.

The proper orthogonal decomposition is a technique for reducing the dimensionality of a given dynamical
system, by projecting it onto a space spanned by the orthonormal columns of a matrix ${\sf V}_k$. To this end, we consider a set of \emph{snapshot solutions} $\bu_j = \bu(t_j)$ of the system \cref{vecode} at $n_s$ different time 
instances ($0 \leq t_1 < \cdots < t_{n_s} \leq T_f$). Let
\begin{equation} \label{linsvd}
\mathsf{S} = [\bu_1, \cdots, \bu_{n_s}] \in \mathbb{R}^{N \times n_s},
\end{equation}
and $\mathcal{S} = \mbox{\tt Range}({\sf S})$ of dimension $d$.
A POD basis of dimension $k < d$ is a set of orthonormal vectors whose 
linear span gives the best approximation, according to some criterion, of the space ${\cal S}$. 
In the 2-norm, this basis can be obtained 
through the singular value decomposition (SVD) of the matrix ${\sf S}$, ${\sf S} = {\sf V}{\sf {\bm \Sigma}}{\sf W}^{\top}$,
with ${\sf V}$ and ${\sf W}$ orthogonal matrices and ${\sf {\bm \Sigma}}={\rm diag}(\sigma_1, \ldots, \sigma_{n_s})$ 
diagonal with non-increasing positive diagonal elements.
If the diagonal elements of ${\sf {\bm \Sigma}}$ have a rapid decay, the first $k$ columns of ${\sf V}$ (left singular vectors)
  are the most dominant in the approximation of ${\sf S}$. Denoting with ${\sf S}_k = {\sf V}_k{\sf {\bm \Sigma}}_k{\sf W}_k^{\top}$
the reduced SVD where only the $k\times k$ top left portion of ${\sf {\bm \Sigma}}$ is retained and ${\sf V}, {\sf W}$ are
truncated accordingly, then $\|{\sf S} - {\sf S}_k\| = \sigma_{k+1}$ \cite{golub13}.

Once the matrix ${\sf V}_k$ is obtained, for $t \in [0, T_f]$ the vector $\bu (t)$ is approximated as 
$\bu (t) \approx {\sf V}_k{\bm y}_k(t)$, 
where the vector ${\bm y}_k(t) \in \RR^k$ solves the {\em reduced} problem
 \begin{equation}
 \dot{{\bm y}}_k(t) = {\bm L}_k{\bm y}_k(t) + {\bm f}_k({\bm y}_k,t) , \quad {\bm y}_k(0) = {\sf V}_k^{\top}\bu_{0}.
 \label{vecodesmall}
 \end{equation}
 Here ${\bm L}_k = {\sf V}_k^{\top}{\bm L}{\sf V}_k$ and ${\bm f}_k({\bm y}_k,t) = 
{\sf V}_k^{\top}{\bm f}({\sf V}_k{\bm y}_k,t)$. 
Although for $k \ll N$ problem \cref{vecodesmall} is cheaper to solve than the original one,
the definition of ${\bm f}_k$ above requires the evaluation of ${\bm f}({\sf V}_k{\bm y}_k,t)$ at each
timestep and
at all $N$ entries, thus still depending on the original system size.
One way to overcome this problem is by means of DEIM.
 
The DEIM procedure, originally introduced in \cite{chaturantabut2010nonlinear}, is utilized to approximate a 
nonlinear vector function ${\bm f}: {\cal T} \rightarrow \RR^N$
by interpolating it onto an empirical basis, that is,
${\bm f}(t) \approx \In {\color{black}{\bm c}}(t),$
where $\{{\bm \varphi}_1,\dots,{\bm\varphi}_p\} \subset \RR^{N}$ is a low dimensional basis,
$\In = [{\bm\varphi}_{1}, \dots, {\bm\varphi}_p] \in \mathbb{R}^{N \times p}$  
and ${\color{black}{\bm c}}(t) \in \mathbb{R}^{p}$ is the vector of time-dependent coefficients to be determined.

Let
$\mathsf{P} = [{\color{black} \bm e}_{\rho_{1}}, \dots, {\color{black} \bm e}_{\rho_{p}}] \in \mathbb{R}^{N \times p}$ be a subset of columns of the identity
matrix, named the ``selection matrix''.
If $\mathsf{P}^{\top}\In$ is invertible, 
in \cite{chaturantabut2010nonlinear} the coefficient vector ${\color{black}{\bm c}}(t)$ is uniquely determined by solving
the linear system
$\mathsf{P}^{\top}\In {\color{black}{\bm c}}(t) = \mathsf{P}^{\top}{\bm f}(t)$, 
so that
\begin{equation} \label{deimsetup}
{\bm f}(t) \approx \In {\color{black}{\bm c}}(t) = \In(\mathsf{P}^{\top}\In)^{-1}\mathsf{P}^{\top}{\bm f}(t).
\end{equation}
The nonlinear term in the reduced model \cref{vecodesmall} is then approximated by
\begin{equation}
\label{fkapprox}
{\bm f}_k({\bm y}_k,t) \approx {\sf V}_k^{\top}\In(\mathsf{P}^{\top}\In)^{-1}\mathsf{P}^{\top}{\bm f}({\sf V}_k{\bm y}_k,t).
\end{equation}

The accuracy of DEIM depends greatly on the basis choice, and in a lesser way by the choice of $\sf P$. 
In most applications the interpolation 
basis $\{{\bm\varphi}_1,\dots,{\bm\varphi}_p\}$ is selected as the POD basis of the set of snapshots 
$\{{\bm f}(\bu_1,t_1), \dots, {\bm f}(\bu_{n_s},t_{n_s})\}$, as described earlier in this section,
that is
given the matrix 
\begin{equation}
\mathsf{N} = [{\bm f}(\bu_1,t_1), \dots, {\bm f}(\bu_{n_s},t_{n_s})] \in \mathbb{R}^{N \times n_s},
\label{nonsvd}
\end{equation} 
the columns of the matrix $\In = [{\bm\varphi}_{1}, \dots, {\bm\varphi}_p]$ are determined as the first $p\le n_s$ 
dominant left singular vectors in the SVD of $\mathsf{N}$.  
The matrix $\sf P$ for DEIM is selected by a greedy algorithm based on the system residual;
see \cite[Algorithm 3.1]{chaturantabut2010nonlinear}. 
In \cite{gugercin2018} the authors showed that a pivoted 
QR-factorization of $\In^{\top}$ may lead to better accuracy and stability properties of
the computed matrix ${\sf P}$. {The resulting approach, called Q-DEIM, 
will be used in the sequel and is implemented as algorithm {\tt q-deim} in \cite{gugercin2018}}.

DEIM is particularly advantageous when the function ${\bm f}$ is evaluated componentwise, in which case it holds that 
$\mathsf{P}^{\top}{\bm f}({\sf V}_k{\bm y}_k,t) = {\bm f}({\sf P}^{\top}{\sf V}_k{\bm y}_k,t)$, so that
the nonlinear function is evaluated only at $p\ll N$ components. 
In general, this occurs whenever finite differences are the underlying discretization method.
If more versatile discretizations 
are required, then different procedures need to be devised.
A discussion of these approaches is deferred to section~\ref{sec:others}.

 \section{A new two-sided proper orthogonal decomposition}\label{sec:matpod}
We derive a \podeim\ algorithm that fully lives in the matrix setting, without
 requiring a mapping from $\RR^{n_x \times n_y}$ to $\RR^N$, so that
no vectors of length $N$ need to be processed or stored.
%
We determine the left and right reduced 
space bases that approximate the space of the given snapshots ${\bm  \Xi}(t)$ (either the nonlinear functions
or the approximate solutions), 
{so that
${\bm  \Xi}(t) \approx {\bm V}_\ell {\bm \Theta}(t)  {\bm W}_{r}^{\top}$, 
where ${\bm V}_\ell, \bm{W}_r$ have $\nu_\ell$ and $\nu_r$ orthonormal columns, respectively.
Their ranges approximate the span of 
the rows (left) and columns (right) spaces of the function ${\bm  \Xi}(t)$, independently
of the time $t$. In practice, we wish to have $\nu_\ell \ll n_x, \nu_r\ll n_y$ so that ${\bm  \Theta}(t)$ will have
a reduced dimension.
%
A simple way to proceed would be to collect all snapshot matrices in a single large (or tall) matrix
and generate the two orthonormal bases corresponding to the rows and columns spaces. This is way too
expensive, both in terms of computational costs and memory requirements.
Instead, we present a two-step procedure 
that avoids explicit computations with all snapshots simultaneously.
The first step sequentially selects the most important information of each
snapshot matrix, relative to the other snapshots, while the second step
prunes the generated spaces by building the corresponding orthonormal bases. 
These two steps can be summarized as follows:

\begin{enumerate}
\item {\it Dynamic selection}. Assume $i$ snapshots have been processed and dominant SVD information
retained. 
For the 
next snapshot ${\bm  \Xi}(t_{i+1})$ perform a reduced SVD and retain the leading singular
triplets in a way that the retained singular values are at least as large as those already kept from previous iterations.
Make sure that at most $\kappa$ SVD components are retained overall, with $\kappa$ selected a-priori;
\item {\it Bases pruning}. Ensure that the vectors spanning the reduced right and left spaces have 
orthonormal columns. Reduce the space dimension if needed.
\end{enumerate}

In the following we provide the details for this two-step procedure.
The strategy that leads to the selection of the actual time instances used
for this construction will be discussed in section~\ref{sec:dynamic}.
To simplify the presentation, and without loss of generality, we assume $n_x=n_y\equiv n$.

\vskip 0.1in
{\it First step.} 
Let $\bm\Xi_i = {\bm \Xi}(t_i)$.
For the collection of snapshots, we wish to determine a (left) reduced basis 
for the range of the matrix $\bm{H}\{{\bm  \Xi} \} := 
({\bm\Xi}_1, \, \ldots, \, {\bm\Xi}_{n_s} ) \in \mathbb{R}^{n \times (n\cdot n_s)}$ and a (right)
reduced basis for the range of the  matrix
$\bm{Z}\{ {\bm  \Xi} \} = ( {\bm\Xi}_1^\top, \ldots, {\bm\Xi}_{n_s}^\top)^\top
= ( {\bm\Xi}_1; \ldots; {\bm\Xi}_{n_s})$. 
This is performed by incrementally including leading components of 
the snapshots ${\bm\Xi}_i$, so as to determine the approximations
\begin{eqnarray*}
{\bm H}\{\bm  \Xi\} &\approx& \widetilde{\bm H}\{\bm  \Xi\} :=  
\widetilde{\bm  V}_{n_s}\widetilde{\bm  \Sigma}_{n_s}\widetilde{\bm  W}_{n_s}^{\top}, 
\qquad \widetilde{\bm  V}_{n_s} \in \mathbb{R}^{n \times \kappa}, \widetilde{\bm  W}_{n_s} \in \RR^{n\cdot n_s \times \kappa} \\
{\bm Z}\{ {\bm  \Xi} \} &\approx& \widehat{\bm  Z}\{ {\bm  \Xi} \} :=
\widehat{\bm V}_{n_s} \widetilde{\bm \Sigma}_{n_s}\widehat{\bm W}_{n_s}^{\top}, 
\qquad \widehat{\bm  V}_{n_s} \in \mathbb{R}^{n\cdot n_s \times \kappa}, 
\widehat{\bm  W}_{n_s} \in \RR^{n\times \kappa}.
\end{eqnarray*}
A rank reduction of the matrices $\widetilde{\bm  V}_{n_s}$ and $\widehat{\bm  W}_{n_s}$ will provide the
sought after bases, to be used for time instances other than those of the snapshots.
Let $\kappa\le n$ be the chosen maximum admissible dimension for the
reduced left and right spaces.
For $i \in \{1, \ldots, n_s\}$, let 
\begin{eqnarray}\label{eqn:Xi_i}
{\bm  \Xi}_i \approx {\bm V}_i{\bm \Sigma}_i{\bm W}_i^\top, \quad {\bm V}_i, {\bm W}_i \in \RR^{n \times \kappa}, 
\quad
{\bm \Sigma}_i = {\rm diag}({\sigma_1^{(i)}, \ldots, \sigma_{\kappa}^{(i)}})
\end{eqnarray}
be the reduced  SVD of ${\bm  \Xi}_i$ corresponding to its dominant $\kappa$ singular
triplets, with singular values
sorted decreasingly. 
Let $\widetilde {\bm  V}_1={\bm  V}_1$, 
$\widetilde {\bm  \Sigma}_1={\bm  \Sigma}_1$ and $\widetilde {\bm  W}_1={\bm  W}_1$.
 For each subsequent $i=2, \ldots, n_s$, 
the leading 
singular triplets of ${\bm   \Xi}_i$ are {\it appended} to the previous matrices
$\widetilde{\bm  \Sigma}_{i-1}$, $\widetilde{\bm  V}_{i-1}$ and $\widetilde{\bm  W}_{i-1}$,
that is
\begin{equation}\label{htilde1}
\widetilde{ {\bm H}}_{i}\{{ \bm  \Xi} \}  = 
  ( \widetilde{\bm V}_{i-1}, \,  {\bm V}_{i} )
\begin{pmatrix}
    \widetilde{{\bm \Sigma}}_{i-1} &  \\
     & {{\bm \Sigma}}_{i}
\end{pmatrix}
\begin{pmatrix}
    \widetilde{\bm W}_{i-1}^{\top}  & \\
    & {\bm W}_{i}^{\top}
  \end{pmatrix}  \equiv \widetilde{\bm V}_{i} {\widetilde{\bm \Sigma}}_{i} \widetilde{\bm W}_{i}^{\top} .
\end{equation}
After at most $n_s$ iterations we will have the final 
$\widetilde{ {\bm H}}\{{ \bm  \Xi} \}$, with no subscript.
The expansion is performed 
ensuring that the appended singular values are at least as large as those already present, so that
the retained directions are the leading ones among all snapshots.
Then the three matrices $\widetilde{\bm V}_{i}$,
$\widetilde{\bm W}_{i}$ and $\widetilde{\bm\Sigma}_{i}$ are truncated so that
the retained diagonal entries of $\widetilde{\bm \Sigma}_i$ are its $\kappa$ largest diagonal elements.
The diagonal elements of $\widetilde {\bm \Sigma}_{i}$ are not
the singular values of ${\bm H}_{i}\{{ \bm  \Xi}\}$, however the adopted truncation strategy
ensures that the error committed is not larger than
$\sigma_{\kappa+1}$, which we define as
 the largest singular value discarded during the whole accumulation process.
Since each column of $\widetilde {\bm V}_{n_s}$ has unit norm, it follows that 
 $\|\widetilde {\bm V}_{n_s}\| \le \kappa$. Moreover,
the columns of $\widetilde{{\bm W}}_{n_s}$ are orthonormal. Hence
$\|{\bm H}\{{ \bm  \Xi}\} - \widetilde{\bm H}\{{ \bm  \Xi}\}\| \le \kappa\sigma_{\kappa+1}$.
In particular, this procedure is not the
same as taking the leading singular triplets of each snapshot per se: this first step allows us to
retain the leading triplets of each snapshot {\it when compared with
all snapshots}, using the magnitude of all singular values as quality measure.}

A similar strategy is adopted to construct the right basis. Formally,
\begin{equation}  \label{ztilde}
\widehat{ {\bm Z}}_{i}\{{ \bm  \Xi} \}  = 
\begin{pmatrix}
    \widehat{\bm V}_{i-1}  & \\
    & {\bm V}_{i}
 \end{pmatrix}   
\begin{pmatrix}
    \widetilde{{\bm \Sigma}}_{i-1} &  \\
     & {{\bm \Sigma}}_{i}
\end{pmatrix}
  \begin{pmatrix} \widehat{\bm W}_{i-1}^\top \\   {\bm W}_{i}^\top \end{pmatrix} 
=
\widehat{\bm V}_{i} {\widetilde {\bm \Sigma}}_{i} \widehat{\bm W}_{i}^{\top} ;
\end{equation}
notice that ${\widetilde {\bm \Sigma}}_{i}$ is the same for 
both $\widetilde {\bm  H}$  and $\widehat {\bm  Z}$, and that
the large matrices $\widetilde {\bm W}_i$ and $\widehat{\bm V}_i$ are not stored explicitly
in the actual implementation.
At completion, the following two matrices are bases candidates for
the left and right spaces:
\begin{equation}\label{eq:basis}
{\rm Left:} \quad \widetilde{\bm V}_{n_s} {\widetilde {\bm \Sigma}}_{n_s}^{\frac 1 2}, 
\qquad
{\rm Right:}\quad {\widetilde {\bm \Sigma}}_{n_s}^{\frac 1 2} \widehat {\bm W}_{n_s}^{\top} ,
\end{equation}
where the singular value matrices keep track of the relevance of each
collected singular vector, and the square root allows us to maintain the
order of magnitude of the snapshot matrices, when the product of the
two left and right matrices is carried out.
Here $n_s$ is the total number of snapshots included in the whole procedure, using the dynamic
procedure discussed in section~\ref{sec:dynamic}.

The procedure is described in Algoritm~\ref{alg:snapstep}.

\begin{algorithm}
\caption{\sc Dynamic selection procedure}
\label{alg:snapstep}
\begin{algorithmic}[1]
\STATE{\textbf{INPUT:} ${\bm \Xi}_i$, $\widetilde{\bm V}_{i-1} \in \RR^{n \times \kappa}$, $\widetilde{\bm \Sigma}_{i-1} \in \RR^{\kappa \times \kappa}$, ${\widehat{\bm W}}_{i-1} \in \RR^{n \times \kappa}$, $\kappa$}
\STATE{\textbf{OUTPUT:}
$\widetilde{\bm V}_{i} \in \RR^{n \times \kappa}$, $\widetilde{\bm \Sigma}_{i} \in \RR^{\kappa \times \kappa}$, ${\widehat{\bm W}}_{i} \in \RR^{n \times \kappa}$.}
\STATE{ Compute
$[{\bm V}_i, {\bm \Sigma}_i, {\bm W}_i] = \mbox{\tt svds}\left({\bm \Xi}_i, \kappa \right)$;}
\STATE{Append $\widetilde{\bm V}_{i} \leftarrow (\widetilde{\bm V}_{i-1}, {\bm V}_i)$,
${\widehat{\bm W}}_{i} \leftarrow ({\widehat{\bm W}}_{i-1},{\bm W}_i)$,
$\widetilde{\bm \Sigma}_{i} \leftarrow \mbox{\tt blkdiag} (\widetilde{\bm \Sigma}_{i-1}, {\bm \Sigma}_i)$;}
\STATE{Decreasingly order the entries of (diagonal) $\widetilde{\bm \Sigma}_{i}$ and keep the first $\kappa$;}
\STATE{Order $\widetilde{\bm V}_{i}$ and ${\widehat{\bm W}}_{i}$ accordingly and keep the first $\kappa$ vectors of each;}
\end{algorithmic}
\end{algorithm}

\vskip 0.1in
{\it Second step.} We complete the two-sided approximation of the snapshot functions 
by pruning the two orthonormal bases associated with the representation
(\ref{eq:basis}). Let
\begin{equation} \label{eqn:Wtilde}
\widetilde{\bm V}_{n_s}{\widetilde{\bm \Sigma}_{n_s}}^{\frac 1 2} =
 \overline{{\bm V}}\,\overline{{\bm \Sigma}}
\overline{\bm W}^\top \quad {\rm and} \quad
{\widetilde{\bm \Sigma}_{n_s}}^{\frac 1 2}\widehat{\bm  W}_{n_s}^{\top}
=     
\breve{\bm V}
\breve{\bm \Sigma} \breve{\bm W}^{\top} 
\end{equation}
be the singular value decompositions of the given matrices.
If the matrices $\overline{\bm \Sigma}$ and $\breve{\bm \Sigma}$ 
have rapidly decaying singular values, 
we can further reduce the low rank approximation of each ${\bm \Xi}_i$.
More precisely, let 
\begin{equation}\label{eq:partitioning}
\overline{\bm V}\, \overline{{\bm \Sigma}} = [{\bm V}_\ell, {\bm V}_{\cal E}] \begin{pmatrix}
    \overline{{\bm \Sigma}}_{\ell} &  \\
    & \overline{{\bm \Sigma}}_{{\cal E}}\end{pmatrix}, \quad
{\rm and} \quad
\breve{{\bm \Sigma}} \breve{\bm W}^\top = \begin{pmatrix}
    \breve{\bm \Sigma}_r &  \\
    & \breve{\bm \Sigma}_{\cal E}
\end{pmatrix}
\begin{bmatrix}{\bm W}_r^\top \\ {\bm W}_{\cal E}^\top \end{bmatrix},
\end{equation}
with ${\bm{V}}_{\ell} \in \RR^{n \times \nu_{\ell}}, {\bm W}_r \in \RR^{n \times \nu_{r}}$.
The final reduced dimensions $\nu_\ell$ and $\nu_r$, that is the number of columns to be retained 
in the matrices ${{\bm V}}_{\ell}$ and ${\bm W}_r$, respectively,
 is determined by the following criterion:
\begin{equation} \label{kselect}
\|\overline{\bm \Sigma}_{\cal E}\|_F \le \frac{\tau}{\sqrt{n_{\max}}} \|\overline{\bm \Sigma}\|_F
\quad \mbox{and} \quad 
\|\breve{\bm \Sigma}_{\cal E}\|_F \le \frac{\tau}{\sqrt{n_{\max}}} \|\breve{\bm \Sigma}\|_F ,
\end{equation}
for some chosen tolerance $\tau\in (0,1)$ and $n_{\max}$ is the maximum number of available snapshots 
of the considered function.

We have assumed so far that at least one singular triplet is retained for all ${ \bm  \Xi}$'s. In practice,
it might happen that for some $i$ none of the singular values is large enough to be retained.
In this case, 
${ \bm  \Xi}_i$ does not contribute to the two-sided basis. 

\begin{remark} \label{cor:symmetry}
If ${\bm  \Xi}_i={\bm  \Xi}(t_i)$ is symmetric for $t_i \in [0, T_f]$,
the reduction process can preserve this structure.
Indeed,
since ${\bm  \Xi}_i$ is symmetric 
it holds that
${ \bm  \Xi}_i={{\bm V}_i}{{\bm \Sigma}}_i{{\bm W}_i}^{\top} = {{\bm V}_i}{{\bm \Sigma}}_i{\bm D}_i{{\bm V}_i}^{\top}$, with ${\bm D}_i$ diagonal of ones and
minus ones. As a consequence, ${\bm W}_r={\bm V}_\ell$.
Positive definiteness can also be preserved, with ${\bm D}_i$ the identity matrix.
\end{remark}

In the following we use the pair $({\bm  V}_{\ell}, {\bm {W}}_{r})$,
hereafter denoted as
{\it two-sided  proper orthogonal decomposition} ({\sc 2s-pod)}, to 
approximate the function ${\color{black} \bm  \Xi}(t)$ for some $t\ne t_i$:
\begin{equation}\label{eqn:approx_Xi}
{\bm  \Xi}(t) \approx {\bm  V}_{\ell,{\bm \Xi}} {\bm \Theta}(t)  {\bm {W}}_{r,{\bm \Xi}}^{\top}
\end{equation}
with $\bm \Theta$ depending on $t$, of reduced dimension, and
${\bm V}_{\ell,{\bm \Xi}}$ and ${\bm W}_{r,{\bm \Xi}}$ 
play the role of ${\bm V}_\ell$ and ${\bm W}_{r}$ respectively.

 \section{Connections to other matrix-based interpolation POD strategies}\label{sec:others}
The approximation discussed in the previous section is not restricted to problems
of the form (\ref{Realproblem}), but rather it can be employed to any POD function approximation
where the snapshot vectors are transformed into matrices, giving rise to a matrix DEIM 
methodology. This class of approximation has been explored in the recent literature,
where different approaches have been discussed, especially in connection with parameter-dependent
problems and Jacobian matrix approximation; see, e.g., \cite{wirtz2014},\cite{carlberg2015}, 
and the thorough discussion in \cite{benner2015review}. 
In the former case, the setting is particularly
appealing whenever the operator has a parameter-based affine function formulation,
while in the Jacobian case the problem is naturally stated in matrix terms, possibly with
a sparse structure \cite{bonomi2017},\cite{sandu2017}.
In the nonaffine case, in \cite{bonomi2017},\cite{negri2015} an affine matrix approximation (MDEIM)
was proposed by writing appropriate (local) sparse representations of the POD basis, as is the case 
for finite element methods. As an alternative in this context, it was
shown in \cite{dedden2012}  that DEIM can be applied locally to functions defined on the 
unassembled finite element mesh (UDEIM); we refer the reader to  \cite{tiso2013} for 
more details and to \cite{antil2014} for a detailed experimental analysis.

In our approach we consider the approximation in (\ref{eqn:approx_Xi}).
If ${\bm  \Theta}(t)$ were diagonal, then this approximation could be thought of as an MDEIM
approach, since then ${\color{black} \bm  \Xi}(t)$ would be approximated by a sum of rank-one matrices with
the time-dependent diagonal elements of ${\bm  \Theta}(t)$ as coefficients.
Instead, in our setting ${\bm  \Theta}(t)$ is far from diagonal, hence our approximation determines
a more general memory saving approximation based
on the low rank representation given by ${\bm V}_{\ell,\Xi}, {{\bm W}}_{r,\Xi}$.

Another crucial novel fact of our approach  is the following.
While methods such as MDEIM aim at creating a linear combination of matrices,
they still rely on the vector DEIM for computing these matrices, thus only detecting
the leading portion of the left range space.
In our construction, the left and right approximation spaces spanned by
${\bm V}_{\ell,\Xi}, {{\bm W}}_{r,\Xi}$, respectively, stem from
a subspace selection of the range spaces of the whole snapshot matrix ${\bm  H}\{{\color{black} \bm  \Xi}\}$ (left space)
and  ${\bm  Z}\{{\color{black} \bm  \Xi}\}$ (right space); here both spaces are {\it matrix} spaces.
In this way, the leading components of both spaces can be captured. In particular,
specific space directions taken by the approximated function during the time evolution
can be more easily tracked; we refer to the previous version of this
manuscript for an experimental illustration \cite[section 8.1]{Kirsten.Simoncini.arxiv2020}.

In light of the discussion above,
our approach might also be interpreted in terms of the ``local basis'' POD framework, see, e.g.,
\cite{amsallem2011}, where the generality of the bases is ensured
by interpolation onto matrix manifolds. For a presentation of this methodology we also
refer the reader to the insightful survey \cite[section 4.2]{benner2015review}.
In this context, the matrices ${\bm V}_\ell, {{\bm W}}_{r}$ may 
represent a new {\it truncated} interpolation 
of the matrices in ${\bm  H}\{{\color{black} \bm  \Xi}\}$ and ${\bm  Z}\{{\color{black} \bm  \Xi}\}$,
in a completely algebraic setting.

\section{A dynamic algorithm for creating the \MPOD\ approximation space} \label{sec:dynamic}
We describe an adaptive procedure for selecting the time instances employed
in the first step of the basis construction of section~\ref{sec:matpod}. This procedure
will be used for the selection of both the solution and the nonlinear function snapshots. The 
dynamic procedure starts with a coarse discretization of the time interval (using one forth of the
available nodes), and then continues with two successive refinements if needed.

Let $n_{\max}$ be the maximum number of available snapshots of the considered function ${\bm \Xi}(t)$, $t \in [t_0, T_f]$
with $t_0=0$.
A first set ${\cal I}_1$ of $n_{\max}/4$ equispaced time instances in $[0, T_f]$ is considered
(symbol `$\star$' in \cref{line}). If needed, a second set ${\cal I}_2$ of $n_{\max}/4$ equispaced time instances
are considered (symbol `$\times$'),
whereas the remaining $n_{\max}/2$ time instances (symbol `$\square$' and set ${\cal I}_3$) are 
considered in the third phase, if needed at all. 

\begin{figure}[htb]
\begin{center}
    \includegraphics[width=.35\textwidth]{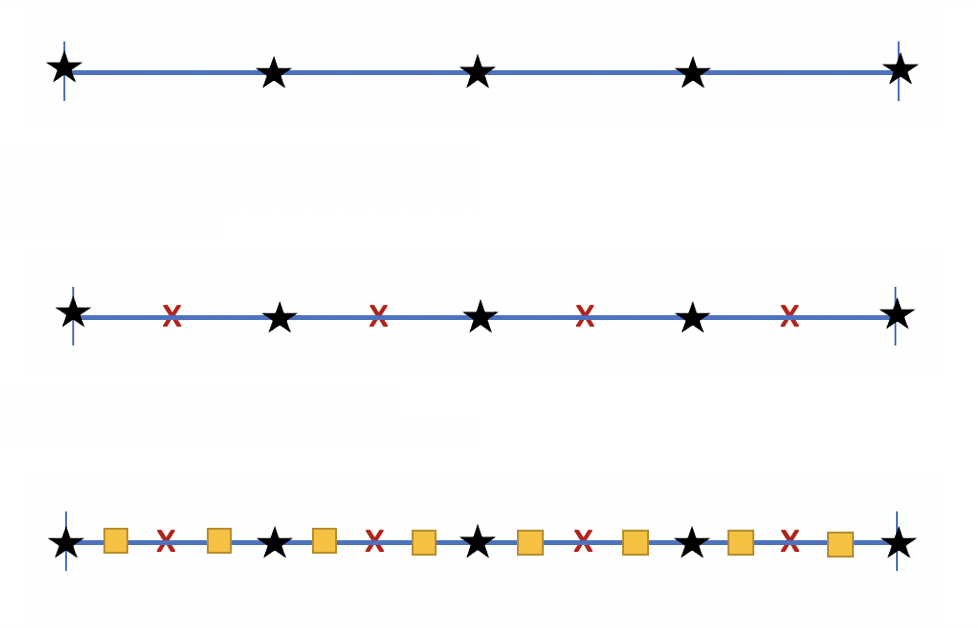}
    \caption{The three evaluation phases of the refinement procedure. \label{line}}
    \end{center}
\end{figure}

The initial \MPOD\ basis matrices of dimension $\kappa$, i.e. $\widetilde{\bm V}_1$ and $\widehat{\bm W}_1$ 
from \cref{eq:basis}, are constructed by processing the snapshot ${\bm \Xi}(t_0)$. 
For all other time instances $t_i$ in each phase, 
we use the following inclusion criterion 
\begin{equation}
{\rm if} \quad \epsilon_i := \frac{\|{\bm \Xi}(t_i) - \Pi_\ell {\bm \Xi}(t_i) \Pi_r\|}{\|{\bm \Xi}(t_i)\|} 
> {\tt tol} \quad {\rm then \, include}
\label{error}
\end{equation}
where $\Pi_\ell$ and $\Pi_r$ are orthogonal projectors onto the left and right spaces (these are 
implicitly constructed
by performing a reduced QR decomposition of the current 
matrices $\widetilde{\bm V}_i$ and $\widehat{\bm W}_i$ on the fly).
If a snapshot is selected for inclusion, the leading
singular triplets of ${\bm \Xi}(t_i)$ are appended to the current
bases, and the leading $\kappa$ components are retained as in Algorithm~\ref{alg:snapstep};
then the next time instance in the phase is investigated. 
 If, by the end of the phase, the {arithmetic}
mean of the errors $\epsilon_i$ in \cref{error} is above  {{\tt tol}}, it means that the bases are not sufficiently good
and we move on to the next refinement phase. Otherwise, the snapshot selection procedure is ended and the matrices $\overline{\bm V}$ and $\breve{\bm W}$ are computed and
pruned by the second step in section~\ref{sec:matpod} to form 
the final {\sc 2s-pod} basis matrices ${\bm V}_{\ell,{\bm \Xi}} \in \RR^{n \times \nu_\ell}$ and ${\bm W}_{r,{\bm \Xi}} \in \RR^{n \times \nu_r}$. 
The full {\sc dynamic} procedure
to create the \MPOD\ approximation space is presented in \cref{alg:snapadap}.}

\begin{algorithm}
\caption{\texttt{{\sc dynamic} \MPOD\ }}
\label{alg:snapadap}
\begin{algorithmic}[1]
\STATE{\textbf{INPUT:}
Function ${\bm \Xi}\,:\, {\cal T} \mapsto \RR^{n \times n}$, $n_{\max}$, {\tt tol}, $\kappa$, phase sets ${\cal I}_{1,2,3}$}
\STATE{\textbf{OUTPUT:}  ${\bm V}_{\ell,{\bm \Xi}} \in \RR^{n \times \nu_{\ell}}$ and ${\bm W}_{r,{\bm \Xi}} \in \RR^{n \times \nu_{r}}$}
\STATE{ Compute $[{\bm V}_1, {\bm \Sigma}_1, {\bm W}_1] = \mbox{\tt svds}\left({\bm \Xi}(t_0), \kappa\right)$;}
\STATE{ Let $\widetilde{\bm V}_1 = {\bm V}_1$, 
$\widehat{\bm W}_1 = {\bm W}_1, \widetilde{\bm \Sigma}_1 = {\bm \Sigma}_1$;}
\vskip 0.05in
\STATE{\hspace{0.1cm} {\it Dynamic selection step :}}
\vskip 0.05in
\FOR{{\sc phase} $= 1,2,3$}
\FOR{all $t_i \in {\cal I}_{\footnotesize\mbox{\sc phase}}$}
\IF{\cref{error} satisfied}
\STATE{Process snapshot ${\bm \Xi}(t_i)$ using Dynamic selection (\cref{alg:snapstep});}
\STATE{Update {$\widetilde{\bm V}_i$ and $\widehat{\bm W}_i$};} 
\ENDIF
\ENDFOR
\IF{$ \sum_{t_i\in {\cal I}_{\tiny\mbox{\sc phase}}} \epsilon_i  \le {\tt tol} |{\cal I}_{\tiny\mbox{\sc phase}}|$}
\STATE{{\bf break} and go to 17;} 
\ENDIF
\ENDFOR
\vskip 0.05in
\STATE{\hspace{0.1cm} {\it Bases Pruning:}}
\vskip 0.05in
\STATE{Determine the reduced SVD of $\widetilde{\bm V}_{n_s}\widetilde{\bm\Sigma}_{n_s}^{\frac 1 2}$
and $\widetilde{\bm\Sigma}_{n_s}^{\frac 1 2} \widehat {\bm W}_{n_s}^\top$ in (\ref{eqn:Wtilde});}
%
\STATE{Determine the final reduced ${\bm V}_{\ell,{\bm \Xi}}$ and ${\bm W}_{r,{\bm \Xi}}$
as in (\ref{eq:partitioning}) using the criterion (\ref{kselect});}
\STATE{\bf Stop} 
\end{algorithmic}
\end{algorithm}

\section{Approximation of the nonlinear function ${\cal F}_k$ in the reduced model} \label{sec:matdeim}
To complete the reduction of the original problem to the small size problem
(\ref{Realproblemsmall0}),
we need to discuss the derivation of the approximation $\reallywidehat{{\cal F}_{k}({\bm Y}_k,t)}$.
Let $\{{\cal F}(t_j)\}_{j=1}^{n_s}$  be a set of snapshots of the nonlinear function ${\cal F}$
at times $t_j$, $j=1, \ldots, n_s$. Using Algorithm~\ref{alg:snapadap} 
 we compute the two matrices ${\bm V}_{\ell,{\cal F}}\in\RR^{n\times p_1}$,
${\bm W}_{r,{\cal F}}\in\RR^{n\times p_2}$ so as to approximate ${\cal F}(t)$ as
\begin{equation}\label{eq:Ftilde}
{\cal F}(t) \approx {\bm V}_{\ell,{\cal F}}{\bm C}(t){\bm W}_{r,{\cal F}}^{\top},
\end{equation}
with ${\bm C}(t)$ to be determined.
Here $p_1, p_2$ play the role of $\nu_{\ell}, \nu_r$ in the general description, and they will be used throughout
as basis truncation parameters for the nonlinear snapshots.
By adapting the DEIM idea to a two-sided perspective,
the coefficient matrix ${\bm C}(t)$ is determined by selecting independent rows from the matrices 
${\bm  V}_{\ell,{\cal F}}$ and ${\bm  W}_{r,{\cal F}}$, so that
$$
{\bm P}_{\ell, {\cal F}}^{\top}{\bm  V}_{\ell,{\cal F}}{\bm C}(t){\bm  W}_{r,{\cal F}}^{\top}{\bm P}_{r,{\cal F}} = {\bm P}_{\ell, {\cal F}}^{\top}{\cal F}(t){\bm P}_{r,{\cal F}} ,
$$
where ${\bm P}_{\ell, {\cal F}} = [e_{\pi_1}, \cdots,e_{\pi_{p_1}}] \in \mathbb{R}^{n \times p_1}$ and 
${\bm P}_{r,{\cal F}} = [e_{\gamma_1}, \cdots,e_{\gamma_{p_2}}]\in \mathbb{R}^{n \times p_2}$ 
are columns of the identity matrix of size $n$.
Both matrices are defined similarly to the selection matrix ${\sf P}$ from section~\ref{sec:POD_DEIM}, 
and they act on ${\bm  V}_{\ell,{\cal F}}, {\bm  W}_{r,{\cal F}}$, respectively.
If $ {\bm P}_{\ell, {\cal F}}^{\top}{\bm  V}_{\ell,{\cal F}}$ and  
${\bm P}_{r,{\cal F}}^{\top}{\bm  W}_{r,{\cal F}}$ are nonsingular, 
then the coefficient matrix ${\bm C}(t)$ is determined by
$$
{\bm C}(t) = ({\bm P}_{\ell, {\cal F}}^{\top}{\bm  V}_{\ell,{\cal F}})^{-1}{\bm P}_{\ell, {\cal F}}^{\top}{\cal F}(t){\bm P}_{r,{\cal F}}({\bm  W}_{r,{\cal F}}^{\top}{\bm P}_{r,{\cal F}})^{-1}.
$$
With this coefficient matrix ${\bm C}(t)$, the final approximation (\ref{eq:Ftilde}) 
becomes\footnote{If the nonlinear function ${\cal F}(t)$ is symmetric for 
all $t \in [0, T_f]$, thanks to \cref{cor:symmetry} this approximation will preserve the symmetry of the 
nonlinear function.}
\begin{equation} \label{deimapprox}
\widetilde{\cal F}(t) =
{\bm  V}_{\ell,{\cal F}}({\bm P}_{\ell, {\cal F}}^{\top}{\bm  V}_{\ell,{\cal F}})^{-1}{\bm P}_{\ell, {\cal F}}^{\top}{\cal F}(t)
{\bm P}_{r,{\cal F}}({\bm  W}_{r,{\cal F}}^{\top}{\bm P}_{r,{\cal F}})^{-1} {\bm  W}_{r,{\cal F}}^{\top} 
=: {\color{black} \bm Q}_{\ell,{\cal F}}{\cal F}(t){\color{black} \bm Q}_{r,{\cal F}}^{\top} .
\end{equation} 
Note that ${\color{black} \bm Q}_{*,{\cal F}}$ are oblique projectors. 
{A similar approximation can be found in \cite{sorensen2016}}.
In addition to that of the two spaces, an important role is played by the
choice of the interpolation indices contained in ${\bm P}_{\ell, {\cal F}}$ and ${\bm P}_{r,{\cal F}}$.
We suggest determining {\color{black}these indices 
for the matrices ${\bm P}_{\ell, {\cal F}}$ and ${\bm P}_{r,{\cal F}}$ as the output of {\tt q-deim} (\cite{gugercin2018})} with inputs ${\bm  V}_{\ell,{\cal F}}$ and ${\bm  W}_{r,{\cal F}}$, respectively.

We next provide a bound measuring the distance between the error obtained with
the proposed oblique projection (\ref{deimapprox}) and the 
best approximation error of ${\cal F}$ in the same range spaces, 
where we recall that ${\bm  V}_{\ell,{\cal F}}$ and ${\bm  W}_{r,{\cal F}}$ have orthonormal columns.
This bound is a direct extension to the matrix setting of \cite[Lemma 3.2]{chaturantabut2010nonlinear}.
\begin{proposition}
\label{errprop}
Let ${\cal F} \in \RR^{n \times n}$ be an arbitrary matrix, and let
$
 \widetilde{\cal F} = {\color{black} \bm Q}_{\ell,{\cal F}}{\cal F}{\color{black} \bm Q}_{r,{\cal F}}^{\top},
$ 
as in \cref{deimapprox}. Then
%
\begin{equation}
\label{errbound}
\| {\cal F} - \widetilde{\cal F} \|_F \le {\color{black} c_{\ell}}{\color{black} c_{r}}\, 
\|{\cal F} - {\bm  V}_{\ell,{\cal F}}{\bm  V}_{\ell,{\cal F}}^{\top} {\cal F}{\bm  W}_{r,{\cal F}}{\bm  W}_{r,{\cal F}}^{\top} \|_F
\end{equation}
where ${\color{black} c_{\ell}} = \left \| ({\bm P}_{\ell, {\cal F}}^{\top}{\bm  V}_{\ell,{\cal F}})^{-1} \right \|_2$ 
and ${\color{black} c_{r}} = \left \| ({\bm P}_{r,{\cal F}}^{\top}{\bm  W}_{r,{\cal F}})^{-1} \right \|_2$.
\end{proposition}

{\it Proof.}
Recall that ${\bm f} = \mbox{\tt vec}({\cal F})$. Then, by the properties of the Kronecker product
{\small \begin{equation*}
\begin{split}
\| {\cal F} - \widetilde{\cal F} \|_F &= \| \mbox{\tt vec}({\cal F}) - \mbox{\tt vec}(\widetilde{\cal F}) \|_2 = \| {\bm f} - ({\color{black} \bm Q}_{r,{\cal F}} \otimes {\color{black} \bm Q}_{\ell,{\cal F}}){\bm f} \|_2 \\ &= \left \|{\bm f} - ({\bm  W}_{r,{\cal F}} \otimes {\bm  V}_{\ell,{\cal F}})\left(({\bm P}_{r,{\cal F}} \otimes {\bm P}_{\ell, {\cal F}})^{\top}({\bm  W}_{r,{\cal F}} \otimes {\bm  V}_{\ell,{\cal F}})\right)^{-1}({\bm P}_{r,{\cal F}} \otimes {\bm P}_{\ell, {\cal F}})^{\top}{\bm f}\right\|_2
\end{split}
\end{equation*}}
Therefore, by \cite[Lemma 3.2]{chaturantabut2010nonlinear}, 
{\footnotesize
\begin{eqnarray*}
{\| {\cal F} - \widetilde{\cal F} \|_F} &\le&
{
     \left\|\left(({\bm P}_{r,{\cal F}} \otimes
{\bm P}_{\ell, {\cal F}})^{\top}({\bm  W}_{r,{\cal F}} \otimes
{\bm  V}_{\ell,{\cal F}})\right)^{-1}\right\|_2
\left\|{\bm f} - ({\bm  W}_{r,{\cal F}} \otimes
{\bm  V}_{\ell,{\cal F}})({\bm  W}_{r,{\cal F}} \otimes {\bm  V}_{\ell,{\cal F}})^{\top}{\bm f}\right\|_2
}
\\
&=& 
{
\left\|({\bm P}_{\ell, {\cal F}}^{\top}{\bm  V}_{\ell,{\cal F}})^{-1}\right\|_2
\left\|({\bm P}_{r,{\cal F}}^{\top}{\bm  W}_{r,{\cal F}})^{-1}\right\|_2 
\left \|{\cal F} - {\bm  V}_{\ell,{\cal F}}{\bm  V}_{\ell,{\cal F}}^{\top} 
{\cal F}{\bm  W}_{r,{\cal F}}{\bm  W}_{r,{\cal F}}^{\top} \right \|_F.}\,\, \square
\end{eqnarray*}
}

We emphasize that ${\color{black} c_{\ell}}, {\color{black} c_{r}}$ do not depend on time, in case $\cal F$ does.
As has been discussed in \cite{chaturantabut2010nonlinear},\cite{gugercin2018}, it is clear from \cref{errbound} that 
minimizing these amplification factors will minimize the error norm with respect to the best approximation onto the
spaces $\mbox{\tt Range}({\bm  V}_{\ell,{\cal F}})$ and $
\mbox{\tt Range}({\bm  W}_{r,{\cal F}})$. The quantities ${\color{black} c_{\ell}}$, ${\color{black} c_{r}}$ depend 
on the interpolation indices. 
If the indices are selected greedily, as in \cite{chaturantabut2010nonlinear}, then
\begin{equation}\label{eqn:bound_P}
{\color{black} c_{\ell}}
\le \frac{(1 + \sqrt{2n})^{p_1 - 1}}{\|e_1^T {\bm  V}_{\ell,{\cal F}}\|_{\infty}}, \quad 
{\color{black} c_{r}} 
\le \frac{(1 + \sqrt{2n})^{p_2 - 1}}{\|e_1^T {\bm  W}_{r,{\cal F}}\|_{\infty}}.
\end{equation}
If the indices are selected by a pivoted 
QR factorization as in {\color{black}\tt q-deim}, then   
\begin{eqnarray*}
  {\color{black} c_{\ell}}\le \sqrt{n-p_1 + 1}\frac{\sqrt{4^{p_1} + 6p_1 -1}}{3}, \quad
    {\color{black} c_{r}} \le \sqrt{n-p_2 + 1}\frac{\sqrt{4^{p_2} + 6p_2 -1}}{3} ,
\end{eqnarray*}
which are better bounds than those in (\ref{eqn:bound_P}), though still rather pessimistic; 
see \cite{gugercin2018}.
 
To complete the efficient derivation of the reduced model in \cref{Realproblemsmall0} we are left with the 
final approximation of ${\cal F}_k$ in (\ref{deimapprox}).
If ${\cal F}$ is evaluated componentwise, as we assume throughout\footnote{For general nonlinear 
functions the theory from \cite[Section 3.5]{chaturantabut2010nonlinear} can be 
extended to both matrices ${\bm  V}_{\ell,U}$ and ${\bm  W}_{r,U}$.},
then ${\bm P}_{\ell, {\cal F}}^{\top} {\cal F}({\bm  V}_{\ell,U}{\bm Y}(t)
{\bm  W}_{r,U}^{\top},t){\bm P}_{r,{\cal F}} =
{\cal F}({\bm P}_{\ell, {\cal F}}^{\top}{\bm  V}_{\ell,U}{\bm Y}(t){\bm  W}_{r,U}^{\top}{\bm P}_{r,{\cal F}},t)$, so that
{\small \begin{align}
\label{nonlinapprox}
{\cal F}_k({\bm Y}_k, t) &\approx
{\bm  V}_{\ell,U}^{\top}{\bm  V}_{\ell,{\cal F}}({\bm P}_{\ell, {\cal F}}^{\top}{\bm  V}_{\ell,{\cal F}})^{-1}{\bm P}_{\ell, {\cal F}}^{\top}
{\cal F}({\bm  V}_{\ell,U}{\bm Y}_k(t)
{\bm  W}_{r,U}^{\top},t){\bm P}_{r,{\cal F}}({\bm  W}_{r,{\cal F}}^{\top}{\bm P}_{r,{\cal F}})^{-1} 
{\bm  W}_{r,{\cal F}}^{\top}{\bm  W}_{r,U} \nonumber\\
& ={\bm  V}_{\ell,U}^{\top}{\bm  V}_{\ell,{\cal F}}({\bm P}_{\ell, {\cal F}}^{\top}{\bm  V}_{\ell,{\cal F}})^{-1}{\cal F}
({\bm P}_{\ell, {\cal F}}^{\top}{\bm  V}_{\ell,U}{\bm Y}_k(t){\bm  W}_{r,U}^{\top}{\bm P}_{r,{\cal F}},t)
({\bm  W}_{r,{\cal F}}^{\top}{\bm P}_{r,{\cal F}})^{-1}{\bm  W}_{r,{\cal F}}^{\top}{\bm  W}_{r,U}\nonumber \\
& =: \reallywidehat{{\cal F}_k({\bm Y}_k, t)}. 
\end{align}
}

The matrices ${{\bm  V}_{\ell,U}^{\top}{\bm  V}_{\ell,{\cal F}}({\bm P}_{\ell, {\cal F}}^{\top}
{\bm  V}_{\ell,{\cal F}})^{-1}} \in \RR^{k_1 \times p_1}$ 
and ${({\bm  W}_{r,{\cal F}}^{\top}{\bm P}_{r,{\cal F}})^{-1} {\bm  W}_{r,{\cal F}}^{\top}
{\bm  W}_{r,U}} \in \RR^{p_2 \times k_2}$ 
are independent of $t$, therefore they can be precomputed and stored once for all.
Similarly for the products ${\bm P}_{\ell, {\cal F}}^{\top}{\bm  V}_{\ell,U} \in \RR^{p_1 \times k_1}$ 
and ${\bm  W}_{r,U}^{\top}{\bm P}_{r,{\cal F}} \in \RR^{k_2 \times p_2}$. Note that products
involving the selection matrices ${\bm P}$'s are not explicitly carried out: the operation simply requires selecting corresponding rows or columns in the 
other matrix factor.

Finally, we remark that in some cases the full space approximation 
matrix may not be involved. For instance, if
${\cal F}$ is a matrix function (\cite{higham2008}) and $\BU(t)$ is symmetric for all $t \in [0,T_f]$,
{\color{black}so that $\BU(t) \approx {\bm  V}_{\ell,U}{\bm Y}_k(t) {\bm  V}_{\ell,U}^{\top}$, 
then, recalling \cref{eqn:F}, it holds that
$$
{\cal F}_{k}({\bm Y}_k,t) = 
{\bm  V}_{\ell,U}^{\top}{\cal F}({\bm  V}_{\ell,U}{\bm Y}_k{\bm  V}_{\ell,U}^{\top},t){\bm  V}_{\ell,U} \overset{\star}{=} {\cal F}({\bm  V}_{\ell,U}^{\top}{\bm  V}_{\ell,U}{\bm Y}_k,t) = {\cal F}({\bm Y}_k,t),
$$
where the equality $\overset{\star}{=}$ is due to \cite[Corollary 1.34]{higham2008}.}

 \section{Two-sided POD-DEIM for nonlinear matrix-valued ODEs} \label{extend}
To complete the derivation of the numerical method, we need to determine the 
time-dependent matrix ${\bm Y}_k(t)$, $t \in [0,T_f]$ in the approximation 
${\bm  V}_{\ell,U}{\bm Y}_k(t){\bm  W}_{r,U}^{\top} \approx \BU(t)$,
where ${\bm  V}_{\ell,U}\in \mathbb{R}^{n \times k_1}$ and ${\bm  W}_{r,U}\in \mathbb{R}^{n \times k_2}$,
$k_1, k_2 \ll n$ and we let $k=(k_1,k_2)$.
The function ${\bm Y}_k(t)$ is computed as the numerical solution to the reduced problem
(\ref{Realproblemsmall0}) with $\reallywidehat{{\cal F}_k({\bm Y}_k, t)}$ defined in (\ref{nonlinapprox}).
{To integrate the reduced order model \cref{Realproblemsmall0} as time $t$ varies, several
alternatives can be considered. 
The vectorized form of the semilinear problem can be treated with classical first or second
order semi-implicit
methods such as IMEX methods (see, e.g.,\cite{ascher1995}), that appropriately handle the stiff and non-stiff parts of the 
equation.
Several of these methods were originally constructed as rational approximations to exponential integrators,
which have for long time been regarded as too expensive for practical purposes. Recent advances
in numerical linear algebra have allowed a renewed interest in these powerful methods, see, e.g., 
\cite{brachet2020comparaison, garcia2014, Grooms2011}. 
One of the advantages of our matrix setting is that exponential integrators can be
far more cheaply applied than in the vector case, thus allowing for a better treatment 
of the stiff component in the solution; see, 
e.g., \cite{Autilia2019matri}.
Let $\{{\mathfrak t}_i\}_{i=0, \ldots, n_{\mathfrak t}}$ be the nodes discretizing the time interval $[0, T_f]$
with meshsize $h$.
Given the (vector) differential equation $\dot\by={\BL}\by + f(t,\by)$, 
the (first order)
exponential time differencing (ETD) Euler method is given by the following recurrence,
$$
\by^{(i)} = e^{h {\BL}} \by^{(i-1)} + h \varphi_1(h {\BL}) f({\mathfrak t}_{i-1},\by^{(i-1)}),
$$
where $\varphi_1(z) = (e^z -1)/z$. In our setting, ${\BL} = {\bm B}_{k}^\top \otimes \BI_{k_1} +
\BI_{k_2}\otimes {\bm A}_{k}$, for which it holds that
$e^{h {\BL}} = e^{h {\bm B}_{k}^\top}\otimes e^{h {\bm A}_{k}}$ \cite[Th.10.9]{higham2008}, so that
the computation of the exponential of the large matrix
${\BL}$ reduces to
$e^{h{\BL}} {\tt vec}(\BY_k) = {\tt vec}(e^{h {\bm A}_k} \BY_k e^{h {\bm B}_k})$. 
The two matrices used in the exponential functions have now small dimensions,
so that the computation of the matrix exponential is fully affordable. 
As a consequence, the integration step can be performed all at the matrix level, 
without resorting to the vectorized
form. More precisely (see also \cite{Autilia2019matri}), let 
${\mathfrak f}({\bm Y}_k^{(i-1)}) =  \reallywidehat{{\cal F}_k({\bm Y}_k^{(i-1)}, {\mathfrak t}_{i-1})}$.
Then we can compute $\BY_k^{(i)}\approx \BY_k({\mathfrak t}_{i})$ as
$$
\BY_k^{(i)} = e^{h {\bm A}_{k}} \BY_k^{(i-1)} e^{h {\bm B}_{k}} + h\BPhi^{(i-1)},
$$
where the matrix $\BPhi^{(i-1)}$ solves the following linear (Sylvester) matrix equation 
\begin{equation}\label{eq:Sylvester}
{\bm A}_{k} \BPhi + \BPhi {\bm B}_{k} 
= e^{h{\bm A}_{k}} {\mathfrak f}({\bm Y}_k^{(i-1)}) e^{h {\bm B}_{k}} - {\mathfrak f}({\bm Y}_k^{(i-1)}) .
\end{equation}
At each iteration, the application of the matrix exponentials in ${\bm A}_{k}$ and ${\bm B}_{k}$ is
required, together with the solution of the small dimensional Sylvester equation. This linear
equation has a unique solution if and only if the spectra of ${\bm A}_{k}$ and $-{\bm B}_{k}$
are disjoint, a hypothesis that is satisfied in our setting. The solution $\BPhi^{(i-1)}$ can be obtained
by using the Bartels-Stewart method \cite{bartels1972}. In case of symmetric ${\bm A}_{k}$ and ${\bm B}_{k}$,
significant computational savings can be obtained by
computing the eigendecomposition of these two matrices at the beginning of the
online phase, and then compute both the exponential and the solution to the Sylvester equation
explicitly; we refer the reader to \cite{Autilia2019matri} for these implementation details.
Matrix-oriented first order IMEX schemes could also be employed (see \cite{Autilia2019matri}), however 
in our experiments ETD provided smaller errors at comparable computational costs.

{Concerning the quality of our approximation,
error estimates for the full \podeim\ approximation of systems
of the form \cref{vecode} have been derived in \cite{wirtz2014, chat2012}, which also
take into account the error incurred in the numerical solution of the reduced problem.
%
A crucial hypothesis in the available literature is that 
${\bm f}(\bu, t) = \mbox{\tt vec}({\cal F}(\BU,t))$ be Lipschitz continuous with respect to the first argument,
and this is also required for exponential integrators. 
This condition is satisfied for the nonlinear function of our reduced problem. Indeed,
consider the vectorized approximation space 
$\mathbb{V}_{\color{black} \bu} = {\bm W}_{r,{\BU}} \otimes {\bm V}_{\ell,{\BU}}$ and the 
oblique projector  $\mathbb{ Q}_{\bm f} = {\bm Q}_{r,{\cal F}} \otimes {\bm Q}_{\ell,{\cal F}}$
from (\ref{deimapprox}).
If we denote by $\widehat{\bm Y}_k(t)$ the approximate solution of \cref{Realproblemsmall0} with ${\cal F}_k$ approximated
above, then
\begin{equation}
\label{errorequal}
\|{\bm U}(t) - {\bm V}_{\ell,\BU}\widehat{\bm Y}_k(t){\bm W}_{r,{\BU}}^{\top}\|_F = \|\bu(t) - \mathbb{V}_{\color{black} \bu}\widehat{\bm y}_k(t)\|_2,
\end{equation}
where $\widehat{\bm y}_k(t) = \mbox{\tt vec}(\widehat{\bm Y}_k(t))$ solves the reduced problem
$\dot{\widehat{\bm y}}_k(t) = \mathbb{V}_{\color{black}\bu}^{\top}{\bm L}\mathbb{V}_{\color{black} \bu}\widehat{\bm y}_k(t) 
+ \mathbb{V}_{\color{black} \bu}^{\top}\mathbb{Q}_{\color{black} \bm f}{\bm f}(\mathbb{V}_{\color{black} \bu}\widehat{\bm y}_k, t)$.

For ETD applied to semilinear differential equations
the additional requirement is that $\BL$ be sectorial, which can also be assumed for the discretization of the
considered operator $\ell$; see \cite[Th.4 for order $s=1$]{hochbruck2005}.
The error in \cref{errorequal} can therefore be approximated by applying the 
a-priori \cite{chat2012} or a posteriori \cite{wirtz2014} error estimates to the vectorized system 
associated with $\widehat{\bm y}_k$. 
Moreover, if ${\bm f}(\bu, t)$ is Lipschitz continuous, this property
is preserved by the reduced order vector model, since $\mathbb{V}_{\color{black} \bu}$ has orthonormal columns and $\|\mathbb{Q}_{\color{black} \bm f}\|$ is a bounded constant, as shown in \cref{errprop}; see e.g., \cite{chat2012}.

The complete \MPDEIM\ method for the semilinear differential problem \cref{Realproblem} is presented in
Algorithm \MPDEIM. In Table~\ref{tableparameters} we summarize the key dimensions and parameters
of the whole procedure.  A technical
discussion of the algorithm and its computational complexity is postponed to \cref{sec:alg}.

\vskip 0.1in
 \begin{algorithm}
{ {\bf Algorithm}  \MPDEIM\ } 
\hrule
\vskip 0.1in
{\bf INPUT:} Coefficient matrices of \cref{Realproblem}, ${\cal F}: \RR^{n \times n} \times [0,T_f] \rightarrow \RR^{n \times n}$, 
(or its snapshots), 
$n_{\max}$, $\kappa$, and $\tau$, $n_{\mathfrak t}$,
$\{{\mathfrak t}_i\}_{i=0, \ldots, n_{\mathfrak t}}$.

{\bf OUTPUT:} $\BV_{\ell,\BU}, \BW_{r,\BU}$ and ${\bm Y}_k^{(i)}$, $i=0, \ldots, n_{\mathfrak t}$ to implicitly form
the approximation $\BV_{\ell,\BU} {\bm Y}_k^{(i)} \BW_{r,\BU}^\top \approx \BU({\mathfrak t}_i)$

\vskip 0.1in
\hspace{1.2cm}\emph{Offline:}
\begin{enumerate}
\item Determine $\BV_{\ell,\BU}, \BW_{r,\BU}$ for {$\{\BU\}_{i=1}^{n_{\max}}$} and 
$\BV_{\ell,{\cal F}}, \BW_{r,{\cal F}}$ for $\{{\cal F}\}_{i=1}^{n_{\max}}$ via \cref{alg:snapadap} ({\sc dynamic} \MPOD) using
at most $n_s$ of the $n_{\max}$ time instances
(if not available, this includes approximating the snapshots $\{{\cal F}(t_i)\}_{i=1}^{n_{\max}}$, 
$\{{\BU}(t_i)\}_{i=1}^{n_{\max}}$ as the time interval is spanned);
\item Compute ${\bm Y}_k^{(0)} = {\bm  V}_{\ell,U}^{\top}\BU_{0}{\bm  W}_{r,U}$, 
${\bm A}_{k} = {\bm  V}_{\ell,U}^{\top}{\color{black} \bm A}{\bm  V}_{\ell,U}$,
${\bm B}_{k} = {\bm  W}_{r,U}^{\top} {\color{black} \bm B}{\bm  W}_{r,U};$
\item Determine ${\bm P}_{\ell, {\cal F}}, {\bm P}_{r,{\cal F}}$ using {\color{black}\tt q-deim}({\sc 2s-deim});
\item {Compute ${{\bm  V}_{\ell,U}^{\top}{\bm  V}_{\ell,{\cal F}}({\bm P}_{\ell, {\cal F}}^{\top}
{\bm  V}_{\ell,{\cal F}})^{-1}}$, ${({\bm  W}_{r,{\cal F}}^{\top}{\bm P}_{r,{\cal F}})^{-1} {\bm  W}_{r,{\cal F}}^{\top}
{\bm  W}_{r,U}}$,
${\bm P}_{\ell, {\cal F}}^{\top}{\bm  V}_{\ell,U}$ and
${\bm  W}_{r,U}^{\top}{\bm P}_{r,{\cal F}}$};\\

\emph{Online:}
\item For each $i=1, \ldots, n_{\mathfrak t}$
\begin{itemize}
{\item[(i)] Evaluate $
\reallywidehat{{\cal F}_k({\bm Y}_k^{(i-1)},{\mathfrak t}_{i-1})}$ as in \cref{nonlinapprox}
using the matrices computed above;}
\item[(ii)] Numerically solve the matrix equation (\ref{eq:Sylvester}) and compute
$$
\BY_k^{(i)} = e^{h {\bm A}_{k}} \BY_k^{(i-1)} e^{h {\bm B}_{k}} + h\BPhi^{(i-1)} ;
$$
\end{itemize}

\end{enumerate}
\end{algorithm}

\begin{table}[hbt]
\caption{Summary of leading dimensions and parameters of Algorithm \MPDEIM. \label{tableparameters}}
\centering
\begin{tabular}{|l|l|}
\hline
Par.      &  Description  \\ \hline
{\scriptsize $n_s$ }         & {\scriptsize Employed number of snapshots} \\
{\scriptsize $k$ }         & {\scriptsize Dimension of vector POD subspace}      \\
{\scriptsize $p$}   & {\scriptsize Dimension of vector DEIM approx. space}           \\
{\scriptsize $N$}   & {\scriptsize Length of ${\bm u}(t)$, $N = n^2$.}    \\
{\scriptsize $\kappa$}   & {\scriptsize Dimension of the snapshot space approximation }    \\
{\scriptsize $k_i$ }         & {\scriptsize Dimension of left ($i = 1$) and right ($i = 2$) \MPOD\ subspaces}\\
{\scriptsize $p_i$ }         & {\scriptsize Dimension of left ($i = 1$) and right ($i = 2$) \MODEIM\ subspaces}\\
{\scriptsize $n$}   & {\scriptsize Dimension of square ${\bm U}(t)$, for $n=n_x=n_y$}  \\
\hline
\end{tabular}
\end{table}

\section{Numerical experiments} \label{sec:exp}
{In this section we illustrate the performance of our matrix-oriented \MPDEIM\ integrator.
In section~\ref{sec:expesF} we analyze the quality of the approximation space created by the {\sc dynamic} algorithm on 
three 
nonlinear functions {with different characteristics}. Then in section~\ref{sec:full} we focus on the ODE setting by comparing
the new \MPDEIM\ procedure to the standard \podeim.

\subsection{Approximation of a nonlinear function ${\cal F}$}\label{sec:expesF}
We investigate the effectiveness of the proposed 
{\sc dynamic} \MPOD\ procedure for determining the two-sided approximation space of a nonlinear function.

{
As a reference comparison, we consider the vector form of the DEIM approximation (hereafter {\sc vector})
in section~\ref{sec:POD_DEIM}.

We also include comparisons with a simple two-sided matrix reduction strategy that uses
a sequential evaluation of all available snapshots, together with the updating of the
bases ${\bm V}_{\ell,{\cal F}}$ and ${\bm W}_{r,{\cal F}}$ during the snapshot processing.
%
In particular, if $[{\bm V}_j,{\bm \Sigma}_j, {\bm W}_j] = \mbox{\tt svds}\left({\cal F}_j, \kappa \right)$ is the 
singular value decomposition of ${\cal F}_j$ limited to the leading $\kappa$ singular triplets, then
in this simple approach the basis 
matrices ${\bm V}_{\ell,{\cal F}}$ and ${\bm W}_{r,{\cal F}}$ are directly updated by orthogonal reduction of
the augmented matrices
 \begin{equation}\label{truncsvds}
 \begin{pmatrix} {\bm V}_{\ell,{\cal F}}, {\bm V}_j{\bm \Sigma}_j^{\frac 1 2} \end{pmatrix} 
\in \RR^{n \times {\kappa}_1}\quad \mbox{and} \quad 
\begin{pmatrix} {\bm W}_{r,{\cal F}}, {\bm W}_j{\bm \Sigma}_j^{\frac 1 2} \end{pmatrix} \in \RR^{n \times {\kappa}_2}
 \end{equation} 
respectively, where ${\kappa}_1, {\kappa}_2 \ge \kappa$. 
To make this procedure comparable in terms of memory to \cref{alg:snapstep}, we enforce 
that the final dimension $\nu_j$ of each basis satisfies {$\nu_j \le \kappa$, for $j = \ell , r$.}
We will refer to this as the {\sc vanilla} procedure for adding a snapshot to the approximation 
space; see, for instance, \cite{Oxberryetal.17},\cite{Kirsten.21} for additional details.


\begin{example}\label{ex:1}
{\rm Consider the nonlinear functions $\phi_i  \, : \Omega \times [0,T_f] \, \to \RR$, $\Omega \subset \RR^2$, $i=1,2,3$ defined as 
\begin{equation*} 
\footnotesize
\begin{cases}  
\phi_1(x_1,x_2,t) &= \frac{x_2}{\sqrt{ (x_1+x_2 - t)^2 + (2x_1 - 3t)^2 + 0.01^2 }}, \,\, \Omega=[0,2]\times [0,2], T_f=2,\\
\phi_2(x_1,x_2,t) &= \frac{x_1x_2}{(x_2 t+0.1)^2}+ \frac{2^{(x_1+x_2)}}{\sqrt{ (x_1+x_2 - t)^2 + (x_2^2+x_1^2 - t^2)^2 + 0.01^2 }},
\,\, \Omega=[0,1]\times [0,1.5], T_f=3,\\
\phi_3(x_1,x_2,t) &= \frac{x_1(0.1+t)}{(x_2 t+0.1)^2}+ \frac{t2^{(x_1+x_2)}}{\sqrt{ (x_1+x_2 - t)^2 + (x_2^2+x_1^2 - 3t)^2 + 0.01^2 }},
\,\, \Omega=[0,3]\times [0,3], T_f=5.
\end{cases}
\end{equation*}
Each function is discretized with $n=2000$ nodes in each spatial direction to form three matrix 
valued functions ${\cal F}^{(i)}: [0, T_f] \rightarrow \RR^{n \times n}$, for $i = 1,\ldots,3$, respectively. 
{In the truncation criterion (\ref{kselect}),
for all functions we set $n_{max} = 60$ and
${\tau} = 10^{-3}$}. The function $\phi_1$ shows significant variations at the beginning of
the time window, $\phi_3$ varies more towards the right-hand of the time span, while $\phi_2$ is
somewhere in between.

The approximations obtained with the considered methods
are reported in Table~\ref{tablekappa50-70}  for $\kappa = 50$ and $\kappa = 70$,
with the following information: For the adaptive snapshot selection procedure, we indicate the required number of {\sc phases} and
the final total number $n_s$ of snapshot used, the CPU time to construct the basis vectors (time for \cref{alg:snapstep} 
plus time for the SVDs \cref{truncsvds} or {\sc vanilla}), 
the final dimensions $\nu_{\ell}$ and $\nu_{r}$ and finally, the arithmetic mean of the errors 
$
\|{\cal F}(t_j) - {\bm V}_{\ell,{\cal F}}{\bm V}_{\ell,{\cal F}}^{\top}{\cal F}(t_j){\bm W}_{r,{\cal F}}{\bm W}_{r,{\cal F}}^{\top}\|/\|{\cal F}(t_j)\|
$
over 300 equispaced timesteps $t_j$, for each ${\cal F} = {\cal F}^{(i)}$, $i = 1,2,3$.
For the vector approach, where we used $n_s=\kappa$,
the reported time consists of the CPU time needed to perform the
SVD of the long snapshots, while the error is measured using the vector form corresponding to the
formula above; see the description at the beginning of section~\ref{sec:full}.

\begin{table}[htb]
\centering
{\footnotesize
\begin{tabular}{|c|r|r|r|r|r|r|r|r|r|r|r|}
\hline
& & \multicolumn{4}{c|}{$\kappa=50$} & \multicolumn{4}{c|}{$\kappa=70$} \\ \hline
  & & { phases} &  { time} & $\nu_{\ell}$/$\nu_{r}$ & { error}  & {phases} &  { time} & $\nu_{\ell}$/$\nu_{r}$ & { error} \\
${\bm \Xi}$  &{ alg.} & {($n_s$)} &  sec.&  & & {($n_s$)} &  sec.&  & \\ \hline
\multirow{2}{*}{$f_1$}   & {\sc dynamic}  & 2\,\,(9)&    3.5  &  33/39   & $6\cdot10^{-4}$  & 1(7) &    4.7    &  40/50   &    $3\cdot10^{-4}$   \\
& {\sc vanilla}      & -(60) &     27.6    &  42/50 &    $6\cdot10^{-4}$  & -(60) &    38.8    &  42/60 &    $3\cdot10^{-4}$ \\
& {\sc vector}      & -(50) &     35.9    &  41 &    $1\cdot10^{-3}$   & -(70) &     77.3    &  56 &    $3\cdot10^{-4}$\\
\hline
\multirow{2}{*}{$f_2$}   & {\sc dynamic}  & 3(21) &    8.6        &  45/26   &         $8\cdot10^{-4}$  & 2(10) &    6.1  &  48/30   & $6\cdot10^{-4}$ \\
& {\sc vanilla}         & -(60) &          25.6            &  50/37 &       $4\cdot10^{-4}$   & -(60) &          39.6            &  58/37 &       $2\cdot10^{-4}$  \\
& {\sc vector}      & -(50) &     42.7    &  36 &    $6\cdot10^{-3}$   & -(70) &     91.8    &  47 &    $2\cdot10^{-3}$  \\
\hline
\multirow{2}{*}{$f_3$}  & {\sc dynamic}  & 2(11) &    4.4        &  34/33 &   $1\cdot10^{-3}$  & 1(10) &    5.9        &  39/39 &    $3\cdot10^{-4}$  \\
& {\sc vanilla}         &-(60) &          25.4            &  46/46 &       $2\cdot10^{-4}$   &-(60) &          38.6            &  46/46 &       $2\cdot10^{-4}$   \\
& {\sc vector}      & -(50) &     47.1    &  50 &    $2\cdot10^{-3}$   & -(70) &     92.5    &  64 &    $8\cdot10^{-4}$ \\
\hline
\end{tabular}
\caption{Example~\ref{ex:1}. Performance of {\sc dynamic}, {\sc vanilla} and {\sc vector} algorithms for $n = 2000$. \label{tablekappa50-70}}
}
\end{table}

Between the matrix-oriented procedures, the dynamic procedure outperforms the simplified one, both in
terms of space dimensions $n_s$, $\nu_\ell$ and $\nu_r$, and in terms of CPU time, especially when
not all time selection phases are needed. The error is comparable for the two methods. Not unexpectedly,
increasing $\kappa$ allows one to save on the number of snapshots $n_s$, though a 
slightly larger reduced dimension $\nu_{\ell}/\nu_{r}$ may occur.
The vector method is not competitive for any of the observed parameters, taking into account that vectors
of length $n^2$ need be stored.
%
} $\square$
\end{example}

\subsection{Solution approximation of the full problem}\label{sec:full}
We report on a selection of numerical experiments with the {\sc dynamic} \MPDEIM\ procedure on matrix semilinear differential
equations of the form \cref{Realproblem}. 
Once again, we compare the results with a standard vector procedure that applies standard {\sc pod-deim} to the vectorized 
solution and nonlinear function snapshots. In particular, we apply the (vectorized) adaptive procedure 
from \cref{sec:dynamic} using the 
error $\| {\pmb \xi} - {\sf V}_k{\sf V}_k^{\top}{\pmb \xi}\| /\|{\pmb \xi}\|$ for the
selection criterion,
where ${\pmb \xi}$ is the vectorized snapshot and ${\sf V}_k$ is the existing {\sc pod} basis. 
{Notice that in the vector setting we need to process as many nodes, as the 
final space dimension, which depends on the tolerance $\tau$. This is due to the fact that
no reduction takes place.} 

In all experiments CPU times are in seconds, and
all bases are truncated using the criterion in \cref{kselect}. To illustrate the quality of
the obtained numerical solution, we also report on the evaluation of the following average relative error norm 
\begin{equation}\label{errormeasure}
\bar{\cal E}({\bm U})= \frac{1}{n_\mathfrak{t}}
\sum_{\gamma = 1}^{n_\mathfrak{t}}\frac{\| {\bm U}^{(j)} - \widetilde{\BU}^{(j)}\|_F}{\|{\bm U}^{(j)}\|_F},
\end{equation}
where $\widetilde{\BU}^{(j)}$ represents either the {\sc dynamic} approximation, 
or the matricization of the {\sc vector} approximation and 
${\bm U}^{(j)}$ is determined by means of the exponential Euler method applied to the original
problem.

\cref{tableode1} shows the key numbers for the bases construction for both methods. 
For either  ${\bm U}$ or ${\cal F}$
 we report  the number of {\sc phases}, the number $n_s$ of 
included snapshots, and the final space dimensions after the reduction procedure.
We first describe the considered model problem and then comment on the numbers in \cref{tableode1}
and on the detailed analysis illustrated in \cref{breakdowntable}. We stress that these
may be considered typical benchmark problems for our problem formulation.
We also remark that for some of the experimental settings the problem becomes symmetric, so that
thanks to Remark~\ref{cor:symmetry} the two bases could be taken to be the same with further computational
and memory savings. Nonetheless, we do not exploit this extra feature in the reported experiments.

For all examples the full matrix problem has dimension $n = 1000$, while the selected values 
for $\kappa$, $n_{\mbox{\tiny max}}$ and $n_{\mathfrak t}$ are displayed in 
Table \ref{tableode1} and Table \ref{breakdowntable}, respectively.

\begin{example}{\rm {\it The 2D Allen-Cahn equation \cite{allen1979}.} \label{ex2}
Consider the equation\footnote{Note that the linear term $-u$ is kept in the nonlinear 
part of the operator.}
\begin{equation}
\label{ac}
u_{t} = \epsilon_1 \Delta u -  \frac{1}{\epsilon_2^2}\left(u^3-u\right),  \quad 
\Omega = [a,b]^2, \quad t \in [0,T_f],
\end{equation}
with initial condition $u(x,y,0) = u_0$. 
The first example is referred to as {\sc ac 1}. As in \cite{song2016}, we use
the following problem parameters: $\epsilon_1 = 10^{-2}$, $\epsilon_2 = 1$, 
$a= 0$, $b=2\pi$ and $T_f = 5$. Finally, we set $u_0 = 0.05\sin x\cos y$ and 
impose homogeneous Dirichlet boundary conditions.
The second problem (hereafter {\sc ac 2}) is a mean curvature flow problem 
(\cite{evans1991}) of the form \cref{ac} with data as suggested, e.g., in
\cite[Section 5.2.2]{ju2015}, that is
periodic boundary conditions, $\epsilon_1 = 1$, $a = -0.5$, $b= 0.5$,
$T_f = 0.075$, with $u_0 = \tanh \left(\frac{0.4 - \sqrt{x^2 + y^2}}{\sqrt{2}\epsilon_2}\right)$.
Following \cite[Section 5.2.2]{ju2015} we consider 
$\epsilon_2 \in \{0.01, 0.02, 0.04\}$. 
As $\epsilon_2$ decreases stability reasons enforce the use of finer time
discretizations $n_{\mathfrak t}$, also leading to larger values of $n_{\mbox{\tiny max}}$ and
$\kappa$, as indicated in \cref{tableode1} and Table~\ref{breakdowntable}, where we report our experimental results.
} $\square$
\end{example}

\begin{example}{\rm {\it Reaction-convection-diffusion equation.} \label{ex3}
We consider the following reaction-convection-diffusion (hereafter {\sc rcd}) problem, also
 presented in \cite{caliari2009},
\begin{equation}
\label{rda}
u_{t} = \epsilon_1 \Delta u +  (u_x + u_y) +  u(u-0.5)(1-u),  
\quad \Omega = [0,1]^2, \quad t \in [0,0.3] .
\end{equation}
The initial solution is given by 
$u_0 = 0.3 + 256\left(x(1-x)y(1-y)\right)^2$, while zero Neumann boundary 
conditions are imposed.  In \cref{tableode1} and Table~\ref{breakdowntable}
we present results for $\epsilon_1 \in \{0.5, 0.05\}$.
} $\square$
\end{example}

\begin{table}[htb]
\caption{\cref{ex2}: Performance of {\sc dynamic} and {\sc vector} algorithms for $n = 1000$. \label{tableode1}}
\begin{center}
{\footnotesize
\begin{tabular}{|c|r|c|r|rrr|}
\hline
{\sc pb.} &$n_{\mbox{\tiny max}}$/$\kappa$&${\bm \Xi}$ &{\sc algorithm} & {\sc phases} & {$n_s$} &  $\nu_{\ell}$/$\nu_{r}$  \\ \hline
\multirow{4}{*}{{\sc ac 1}}& \multirow{4}{*}{40/50} &\multirow{2}{*}{${\bm U}$}             & {\sc dynamic}        & 1    &8&     9/2     \\
&& & {\sc vector}         & 2    &9 &     9             \\

&& \multirow{2}{*}{${\cal F}$}               & {\sc dynamic}        &1    & 7&    10/3          \\
& && {\sc vector}         & 2   &9 &      9                \\
\hline
\multirow{4}{*}{\shortstack{{\sc ac 2} \\ $\epsilon_2 = 0.04$}}&\multirow{4}{*}{400/50} &\multirow{2}{*}{${\bm U}$}             & {\sc dynamic}        & 1    &2&     15/15     \\
& & &{\sc vector}         & 2   &25&        25          \\

& &\multirow{2}{*}{${\cal F}$}               & {\sc dynamic}        &1    & 3&    27/27          \\
& && {\sc vector}         & 2    &40 &      40                \\
\hline
\multirow{4}{*}{\shortstack{{\sc ac 2} \\ $\epsilon_2 = 0.02$}}&\multirow{4}{*}{1200/70} &\multirow{2}{*}{${\bm U}$}             & {\sc dynamic}        & 1    &3&     30/30     \\
& & &{\sc vector}         & 1    &28&        28          \\

& &\multirow{2}{*}{${\cal F}$}               & {\sc dynamic}        &1   & 4&    39/39          \\
& && {\sc vector}         & 2    &53&      53                \\
\hline
\multirow{4}{*}{\shortstack{{\sc ac 2} \\ $\epsilon_2 = 0.01$}}&\multirow{4}{*}{5000/150} &\multirow{2}{*}{${\bm U}$}             & {\sc dynamic}        & 1    &3&     62/62     \\
& & &{\sc vector}         & 1   &43&        43         \\

& &\multirow{2}{*}{${\cal F}$}               & {\sc dynamic}        &1    & 5&    73/73         \\
& && {\sc vector}         & 2    &92 &      92               \\
\hline
\multirow{4}{*}{\shortstack{ {\sc rdc} \\$\epsilon_1 = 0.5$} }&\multirow{4}{*}{60/50} &\multirow{2}{*}{${\bm U}$}             & {\sc dynamic}        & 1    &3&     10/10     \\
& & &{\sc vector}         & 1   &7 &       7         \\

& &\multirow{2}{*}{${\cal F}$}               & {\sc dynamic}        &1    & 3&    13/13          \\
& && {\sc vector}         & 2    &11 &      11           \\
\hline
\multirow{4}{*}{\shortstack{ {\sc rdc} \\$\epsilon_1 = 0.05$} }&\multirow{4}{*}{60/50} &\multirow{2}{*}{${\bm U}$}             & {\sc dynamic}        & 1    &4&     14/14     \\
& & &{\sc vector}         & 3   &14 &        14          \\

& &\multirow{2}{*}{${\cal F}$}               & {\sc dynamic}        &1    & 3&    17/17          \\
& && {\sc vector}         & 3   &34&      34               \\
\hline
\end{tabular}
}
\end{center}
\end{table}

\begin{table}[htb]
\caption{Computational time and storage requirements of \MPDEIM\   and standard vector
\podeim. CPU times are in seconds and $n= 1000$.\label{breakdowntable}}
\centering
\begin{tabular}{|c|r|rrc|rc|r|}
\hline
& & \multicolumn{3}{|c|}{\sc offline}& \multicolumn{2}{|c|}{\sc online}& \\
\cline{3-5}\cline{6-7}
      &      & {\sc basis} & {\sc deim} & && &{\sc rel.} \\
{\sc pb.}  & {\small \sc method}    & {\sc time} &{\sc time}& {\sc memory} & {\sc time} ($n_{\mathfrak t}$)& {\sc memory} & {\sc error}\\
\hline
\multirow{2}{*}{\sc ac 1} & \sc dynamic      & 1.8    & 0.001  & $200n$ &0.009 (300)& $24n$& $1 \cdot 10^{-4}$   \\
                            &\sc vector      & 0.6   & {0.228} &$18n^2$&0.010 (300)& $18n^2$& $1\cdot 10^{-4}$  \\ \hline
                            \multirow{2}{*}{\shortstack{{\sc ac 2}\\$0.04$}} & \sc dynamic      & 0.8   & 0.005  & $200n$ &0.010 (300)& $84n$& $3 \cdot 10^{-4}$   \\
                            &\sc vector      & {8.4}   & {3.745} &$65n^2$&0.020 (300)& $65n^2$& $2\cdot 10^{-4}$  \\ \hline
                                                        \multirow{2}{*}{\shortstack{{\sc ac 2}\\$0.02$}} & \sc dynamic      & 1.8   & 0.004  & $280n$ &0.140 (1000)& $138n$& $2 \cdot 10^{-4}$   \\
                            &\sc vector      & 14.56   & {5.273} &$81n^2$&0.120 (1000)& $81n^2$& $3\cdot 10^{-5}$  \\ \hline
                                                        \multirow{2}{*}{\shortstack{{\sc ac 2}\\$0.01$}} & \sc dynamic      & 5.3  & 0.008 & $600n$ &0.820 (2000)& $270n$& $5 \cdot 10^{-4}$   \\
                            &\sc vector      & 46.2   & {13.820} &$135n^2$&0.420 (2000)& $135n^2$& $2\cdot 10^{-4}$  \\ \hline
                            \multirow{2}{*}{\shortstack{{\sc rdc}\\$0.5$}} & \sc dynamic      & 0.8   & 0.001  & $200n$ &0.008 (300)& $46n$& $2 \cdot 10^{-4}$   \\
                            &\sc vector      & 0.6   & {0.277} &$18n^2$&0.010 (300)& $18n^2$& $1\cdot 10^{-4}$  \\ \hline
                                                        \multirow{2}{*}{\shortstack{{\sc rdc}\\$0.05$}} & \sc dynamic      & 0.9    & 0.001  & $200n$ &0.010 (300)& $62n$& $2 \cdot 10^{-4}$   \\
                            &\sc vector      & 4.1   & {2.297} &$48n^2$&0.010 (300)& $48n^2$& $1\cdot 10^{-4}$  \\ \hline
\end{tabular}
\end{table}


{\cref{tableode1} shows that for both the solution and the nonlinear function snapshots, the {\sc dynamic} procedure 
requires merely one 
phase, and it retains snapshots at only a few of the time instances. On the other hand, the vector approach typically
requires two or even three phases to complete the procedure. 
The dimension of the bases is not comparable for the matrix and vector approaches, since these are subspaces
of spaces of significantly different dimensions, namely $\RR^n$ and $\RR^{n^2}$ in the matrix and vector cases, respectively.
Nonetheless, it is clear that the memory requirements are largely in favor of the matrix approach, as
shown in \cref{breakdowntable}, where computational and memory details are reported. In particular, in \cref{breakdowntable}
the offline timings are broken down into two main parts. 
The column {\sc basis time}  collects the cost of the Gram-Schmidt orthogonalization in the {\sc vector} setting, and the cumulative cost 
of \cref{alg:snapstep} (for both the solution and nonlinear function snapshots) for 
the {\sc dynamic} procedure. Column {\sc deim time} reports the time required to determine the interpolation indices 
by {\tt q-deim}. 
The ``online times'' report the cost to simulate both reduced order models at $n_{\mathfrak t}$ timesteps.
The relative approximation error in \cref{errormeasure} is also displayed, together with the total memory requirements
for the offline and online parts.
%
In terms of memory,
for the {\sc vector} setting this includes the storage of all the processed snapshots while for \MPDEIM\ that of $\widetilde{\bm V}_i$ and 
$\widehat{\bm W}_i$ from \cref{alg:snapstep}, for both the solution and 
nonlinear function snapshots. This quantity is always equal to $4\kappa\cdot n$ for the {\sc dynamic} setting.

We point out the large gain in basis construction time for the {\sc dynamic} procedure, mainly related
to the low number of snapshots employed (cf. \cref{tableode1}). 
Furthermore, the {\sc dynamic} procedure enjoys a massive gain in memory requirements, for 
very comparable online time and final average errors.

\begin{figure}[htb!]
\centering
\includegraphics[width=0.5\linewidth]{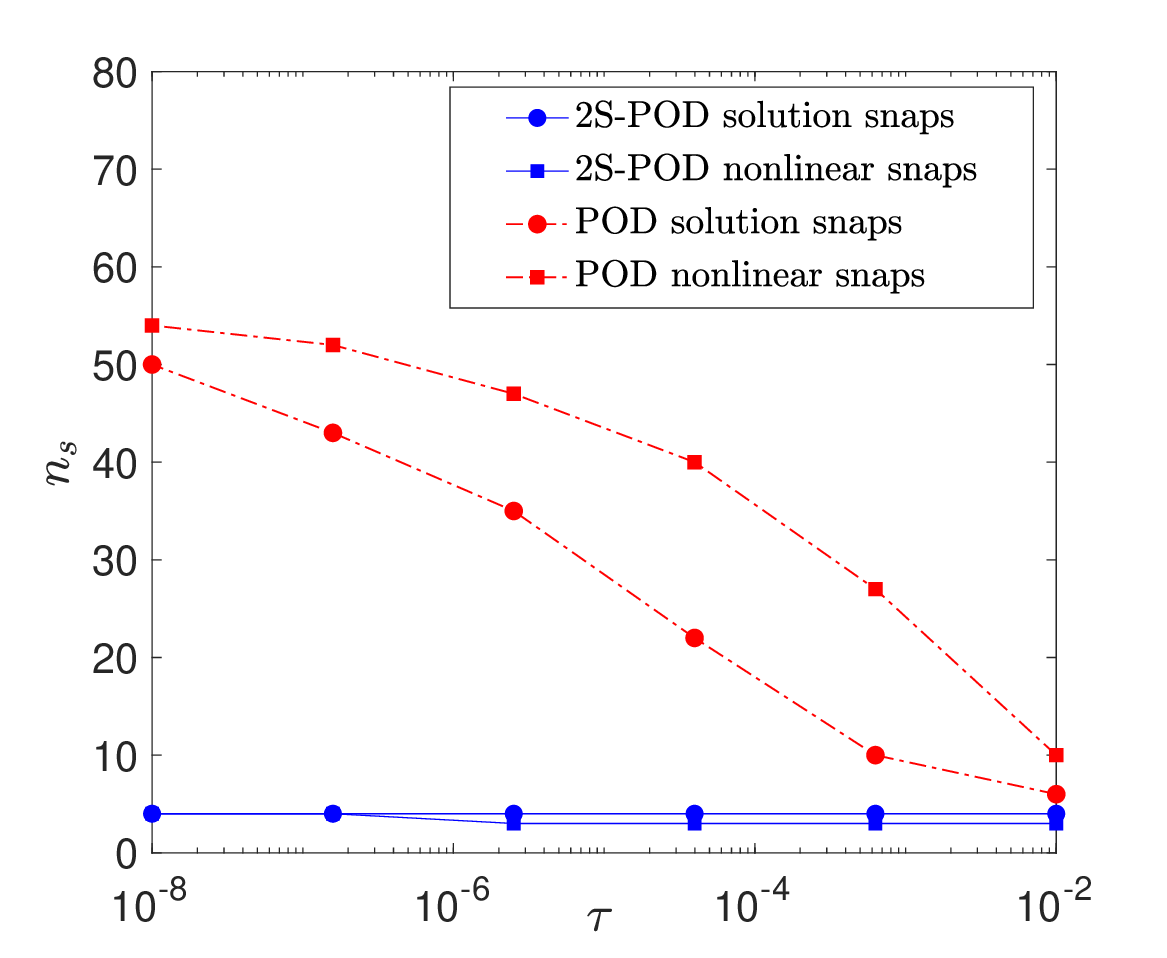}
    \caption{\cref{ex3}. $n = 1000$. 
Number of retained snapshots with respect to $\tau$, for $\epsilon_1=0.05$. \label{snapshotdependence}}
\end{figure}

For the reaction-convection-diffusion example we also
analyze the dependence of the number $n_s$ of retained snapshots on the threshold $\tau$ of the two procedures, having fixed $n_{max} = 60$.
%
For the vector procedure $n_s$
increases as the tolerance $\tau$ decreases, whereas for the 
dynamic procedure $n_s$ remains nearly constant for changing $\tau$. 
This ultimately indicates that the offline cost will increase for the {\sc vector} procedure, if a richer basis is required;
see the appendix and Table~\ref{tablecomplexity} for a detailed discussion.
%
}

 \section{Conclusions and future work} \label{sec:conc}
We have proposed a matrix-oriented \podeim\ type order reduction strategy to 
efficiently handle the numerical solution of
semilinear matrix differential equations in two space variables.
By introducing  a novel interpretation of the proper orthogonal decomposition when 
applied to functions in two variables, we devised a new two-sided discrete interpolation
strategy that is also able to  preserve the symmetric structure in the original nonlinear function
and in the approximate solution. The numerical treatment of the matrix reduced order differential
problem can take full advantage of both the small dimension and the matrix setting, by
exploiting effective exponential integrators.
Our very encouraging numerical experiments show that the new procedure can
dramatically decrease memory and CPU time requirements for the function
reduction procedure in the so-called offline phase. 
Moreover, we illustrated that the reduced low-dimensional matrix differential equation
can be numerically solved in a rapid online phase, without sacrificing too much accuracy. 

This work can be expanded in various directions. In particular,
the companion paper \cite{Kirsten.21} presents a first experimental exploration of the three dimensional
case, which takes advantage of the tensor setting, and uses recently developed tensor linear equations solvers
to advance in time in the (tensor) reduced differential equation; a dynamic approach could
enhance the implementation in \cite{Kirsten.21}.
Generalizations to the multidimensional case and to multiparameters suggest themselves.

\section*{Acknowledgments}
The authors are members of Indam-GNCS, which support is gratefully acknowledged.
Part of this work was also supported by the Grant AlmaIdea 2017-2020 - Universit\`a di Bologna.

\appendix

\section{Discussion of \MPDEIM\ algorithm and computational complexity} \label{sec:alg}
{
{We compare the computational complexity of the 
new \MPDEIM\ method applied to \cref{Realproblem} with that of standard \podeim. All discussions 
are related to Algorithm \MPDEIM\ in \cref{extend}.
\vskip 0.1in

{\it The offline phase.}
The first part of the presented algorithm defines the offline phase. For the $\BU$-set, we considered 
the semi-implicit Euler scheme applied directly to the differential equation
in {\it matrix} form; see, e.g., \cite{Autilia2019matri}, whereas the snapshot selection is done via the adaptive procedure discussed in \cref{sec:dynamic}. 
Notice that moving from one {\sc phase} to the next in \cref{alg:snapadap} does not  require recomputing any quantities. Indeed, if $h^{\ast}$ is the stepsize of the new phase, we determine ${\bm U}(h^{\ast})$ from ${\bm U}(0)$ and initialize the semi-implicit Euler scheme from there; see also \cref{line}.

%

{The computational complexity of approximating the $\kappa$ leading singular triplets with
{\tt svds}, as required by \cref{alg:snapstep}  is given by the implicitly restarted 
Lanczos bidiagonalization, as implemented in the matlab function {\tt svds}. For each $i$ this cost is mainly
given by matrix vector multiplications with the dense $n\times n$ matrix; one Arnoldi cycle involves at most
$2\kappa$ such products, together with $2\kappa$ basis orthogonalizations, leading to ${\cal O}(n^2\kappa + n\kappa)$
operations per cycle \cite{bag2005}.}
{
The final SVDs for {\it bases pruning} at the end of \cref{alg:snapadap} are performed with a dense solver, and  each has 
complexity ${\cal O}(11\kappa^3)$ \cite[p.493]{golub13}. Furthermore, each skinny $QR$-factorization required for the projections ${\Pi_\ell}$ and $\Pi_r$ in \cref{error} has complexity ${\cal O}(2n\kappa^2)$ \cite[p.255]{golub13}. For the standard {\sc pod-deim} algorithm, the reported SVD complexity 
is the total cost of orthogonalizing all selected snapshots by means of Gram-Schmidt, which is ${\cal O}(Nn_s^2)$.

The projected coefficient matrices 
${\bm A}_k$, ${\bm B}_k$, and ${\bm Y}_0$  are computed 
once for all and stored in step 2 of Algorithm \MPDEIM, with a total complexity of 
approximately {${\cal O}(n^2(k_1 + k_2)  +n(k_1 + k_2 + k_1^2 + k_2^2))$}, 
assuming that ${\color{black} \bm A}$ and ${\color{black} \bm B}$ are sparse and ${\BU}_0$ is dense. 
This step is called POD projection in \cref{tablecomplexity}. 

Step 3 in the Algorithm {\sc 2s-pod-deim} has a computational complexity of ${\cal O}(n(p_1^2 + p_2^2))$ \cite{gugercin2018},
{while the matrices in step 4  need to be computed and stored 
with computational complexity {${\cal O}(n(k_1 + k_2)(p_1 + p_2) + (p_1^2 + p_2^2)n + p_1^3 + p_2^3)$}} 
in total \cite{chaturantabut2010nonlinear}. This step is called DEIM projection in \cref{tablecomplexity}.
We recall that the products involving the selection matrices ${\bm  P}_{\ell,{\cal F}}$ and ${\bm  P}_{r,{\cal F}}$
do not entail any computation.

Finally, for the ETD we report the costs for the reduced procedure described in \cite[Section 3.3]{Autilia2019matri}, since we did not experience any stability issues with the diagonalization in our experiments. To this end an a-priori spectral decomposition of each of the reduced matrices 
${\bm A}_k$ and ${\bm B}_k$ is done once for all, 
which has complexity  ${\cal O}(9k_1^3 + 9k_2^3)$ for dense symmetric\footnote{If the coefficient matrices were nonsymmetric, this will be more expensive (still of cubic order), but determining the exact cost is however still an open problem \cite{golub13}.} matrices \cite{golub13}. This makes the cost of evaluating the matrix exponentials negligible, since the computations at each time iteration online will be performed within the eigenbases, following \cite[Section 3.3]{Autilia2019matri}. Furthermore, thanks to the small dimension of the matrices we also explicitly compute the inverse 
of the eigenvector matrices at a cost of ${\cal O}(k_1^3 + k_2^3)$, as required by the online computation. 

All these costs are summarized in \cref{tablecomplexity} and compared with those 
%
of the standard \podeim\ offline phase applied to \cref{vecode}, as indicated in \cite{chaturantabut2010nonlinear}, with dimension $N = n^2$. All coefficient matrices are assumed to be sparse and 
it is assumed that both methods select $n_s$ snapshots via the adaptive procedure. In practice, however, it appears that the two-sided procedure requires far fewer snapshots than the vectorization procedure, as indicated by \cref{tableode1}. 
The table also includes the memory requirements for the snapshots and the basis matrices.

%
\begin{table}[bht]
\caption{Offline phase: Computational costs (flops) for standard 
\podeim\ applied to \cref{vecode},
 and \MPDEIM\ applied to \cref{Realproblem}, { and principal memory requirements. Here $N=n^2$.} \label{tablecomplexity}}
\centering
\begin{tabular}{|l|l|l|}
\hline
Procedure      &  \podeim\ & {\sc dynamic} \MPDEIM\ \\ \hline
{\scriptsize SVD}   & {\scriptsize${\cal O}(Nn_s^2)$  }        & {\scriptsize${\cal O}(n^2\kappa n_s + n\kappa n_s + 6n\kappa^2 + 11\kappa^3)$  }     \\
{\scriptsize QR} & -- & {\scriptsize${\cal O}(n\kappa^2n_s)$  }\\
{\scriptsize DEIM}   & {\scriptsize${\cal O}(Np^2)$}      & {\scriptsize${\cal O}(n(p_1^2 + p_2^2))$ }    \\
{\scriptsize POD projection}  & {\scriptsize${\cal O}(Nk + Nk^2)$ }        & {\scriptsize${\cal O}(n^2(k_1 + k_2)  +n(k_1 + k_2 + k_1^2 + k_2^2))$  }    \\
{\scriptsize DEIM projection}  & {\scriptsize${\cal O}(Nkp + p^2N + p^3 )$}         & {\scriptsize${\cal O}(n(k_1 + k_2)(p_1 + p_2) + (p_1^2 + p_2^2)n + p_1^3 + p_2^3)$}\\
{\scriptsize Snapshot Storage} & {\scriptsize${\cal O}(Nn_s)$} & {\scriptsize${\cal O}(n\kappa)$}    \\
{\scriptsize Basis Storage} & {\scriptsize${\cal O}(N(k + p))$} & {\scriptsize${\cal O}(n(k_1 + k_2 + p_1 + p_2))$}    \\
\hline
\end{tabular}
\end{table}


{\it The online phase.}
The total cost of performing step 5.(i) 
is ${\cal O}(\omega(p_1p_2) + k_1p_1p_2 + k_1k_2p_2)$, where $\omega(p_1p_2)$ 
is the cost of evaluating the nonlinear function at $p_1p_2$ entries. 
Step 5.(ii) requires a matrix--matrix product and the solution of the 
Sylvester equation \cref{eq:Sylvester} in the eigenspace. The latter demands only matrix--matrix products, which come at a cost of ${\cal O}(k_1^2k_2 + k_1k_2^2)$ and two Hadamard products with complexity ${\cal O}(k_1k_2)$; see \cite{Autilia2019matri} for more details. 
This brings the total complexity of one time iteration online to {${\cal O}(\omega(p_1p_2) + k_1p_1p_2 + k_1k_2p_2 + k_1^2k_2 + k_1k_2^2 + k_1k_2)$, which is independent of the original problem size $n$.
}}

}

 \bibliographystyle{siam}
\bibliography{deimbib}
 
\end{document}